\documentclass[a4paper,12pt]{amsart}

\usepackage{times}
\usepackage{amsmath}
\usepackage{amssymb}
\usepackage{latexsym}
\usepackage{color}
\usepackage{fancybox}
\usepackage{ascmac}
\usepackage[pdftex]{graphicx} 
\usepackage{mathrsfs}
\usepackage{amsthm}
\usepackage{enumerate}
\usepackage{cite}

\allowdisplaybreaks[4]

\newcommand{\bbn}{{\mathbb N}}

\newcommand{\bbr}{{\mathbb R}}

\newcommand{\bbz}{{\mathbb Z}}

\newcommand{\al}{{\alpha}}
\newcommand{\be}{{\beta}}

\newcommand{\del}{{\delta}}

\newcommand{\vep}{{\varepsilon}}

\newcommand{\lam}{{\lambda}}
\newcommand{\Lam}{{\Lambda}}

\newcommand{\Om}{{\Omega}}

\newcommand{\cA}{{\mathcal A}}

\newcommand{\cC}{{\mathcal C}}
\newcommand{\cE}{{\mathcal E}}
\newcommand{\cF}{{\mathcal F}}
\newcommand{\cG}{{\mathcal G}}
\newcommand{\cI}{{\mathcal I}}

\newcommand{\cL}{{\mathcal L}}

\newcommand{\cN}{{\mathcal N}}
\newcommand{\cS}{{\mathcal S}}
\newcommand{\cT}{{\mathcal T}}

\newcommand{\cD}{{\mathcal D}}

\newcommand{\cX}{{\mathcal X}}
\newcommand{\cY}{{\mathcal Y}}
\newcommand{\cZ}{{\mathcal Z}}

\newcommand{\sgn}{{\operatorname{sgn}}}

\newcommand{\codim}{\operatorname{codim}}

\newcommand{\D}{\partial}

\theoremstyle{plain}
\newtheorem{thm}{Theorem}[section]
\newtheorem{lem}[thm]{Lemma}
\newtheorem{prop}[thm]{Proposition}

\theoremstyle{definition}

\newtheorem{remk}[thm]{Remark}

\numberwithin{equation}{section}

\makeatletter
\@namedef{subjclassname@2020}{\textup{2020} Mathematics Subject Classification}
\makeatother

\begin{document}
 \pagestyle{plain}

\title{Secondary bifurcations in \\ semilinear ordinary differential equations}

\author{Toru Kan}

\address{Department of Mathematical Sciences, Osaka Prefecture University, 
1-1 Gakuen-cho, Naka-ku, Sakai, 599-8531, Japan
}

\email{kan@ms.osakafu-u.ac.jp}


\begin{abstract}
We consider the Neumann problem for the equation $u_{xx}+\lam f(u)=0$ 
in the punctured interval $(-1,1) \setminus \{0\}$,
where $\lam>0$ is a bifurcation parameter and $f(u)=u-u^3$.
At $x=0$, we impose the conditions 
$u(-0)+au_x(-0)=u(+0)-au_x(+0)$ and $u_x(-0)=u_x(+0)$ for a constant $a>0$
(the symbols $+0$ and $-0$ stand for one-sided limits).
The problem appears as a limiting equation 
for a semilinear elliptic equation in a higher dimensional domain
shrinking to the interval $(-1,1)$.
First we prove that odd solutions and even solutions form families of branches 
$\{ \cC^o_k\}_{k \in \bbn}$ and $\{ \cC^e_k\}_{k \in \bbn}$, respectively.
Both $\cC^o_k$ and $\cC^e_k$ bifurcate from the trivial solution $u=0$.
We then show that $\cC^e_k$ contains no other bifurcation point,
while $\cC^o_k$ contains two points where secondary bifurcations occur.
Finally we determine the Morse index of solutions on the branches.
General conditions on $f(u)$ for the same assertions to hold are also given.
\end{abstract}



\maketitle


\section{Introduction}

This paper is concerned with the problem
\begin{equation}\label{leq}
\left\{
\begin{aligned}
&u_{xx}+\lam f(u)=0, \quad x \in (-1,1) \setminus \{ 0 \}, 
\\
&u_x(-1)=u_x(1)=0, 
\\
&u(-0)+au_x(-0)=u(+0)-au_x(+0),
\\
&u_x(-0)=u_x(+0),
\end{aligned}
\right.
\end{equation}
where $\lam>0$ is a bifurcation parameter, 
$a>0$ is a fixed constant 
and $-0$ (resp. $+0$) stands for the left-hand limit (resp. the right-hand limit) 
as $x$ approaches $0$.
Throughout the paper,
we assume that $f$ satisfies the following conditions:
\begin{equation}\label{basf}
\left\{
\begin{aligned}
&f \in C^2(\bbr), \\
&f(-1)=f(0)=f(1)=0, \ f'(-1)<0, \ f'(0)>0, \ f'(1)<0, \\
&\sgn(u) f(u)>0 \quad \mbox{for} \ u \in (-1,1) \setminus \{0\}, \\
&f(u)=-f(-u) \quad \mbox{for} \ u \in (-1,1).
\end{aligned}
\right.
\end{equation}
Here $\sgn$ is the sign function.
Our interest is the structure of solutions of \eqref{leq} 
in the bifurcation diagram.

\subsection{Background and known results}

The problem \eqref{leq} is related to the stationary problem 
of the scalar reaction-diffusion equation
\begin{equation}\label{rde}
\left\{
\begin{aligned}
&u_t=\Delta u+\lam f(u), && \tilde x \in \Om, \\
&\D_\nu u=0, && \tilde x \in \D \Om,
\end{aligned}
\right.
\end{equation}
where $\Om$ is a bounded domain in a Euclidean space
and $\D_\nu$ stands for the outward normal derivative.
The existence of stable nonconstant stationary solutions
is one of main interests in the study of \eqref{rde}.
In \cite{CH78,Ma79}, 
it was shown that 
\eqref{rde} has no stable nonconstant stationary solutions if $\Om$ is convex,
while in \cite{Ma79}, the existence of such solutions was proved
when $\Om$ is a dumbbell-shaped domain.
Here, by dumbbell-shaped domain, 
we mean a domain given by the union of two regions 
and a thin tubular channel connecting them.
The structure of stable nonconstant stationary solutions was studied
by deriving a finite-dimensional limiting equation
as the thickness of the channel approaches zero.
The limiting equation was obtained and examined by \cite{V83},
and it was shown that a stable nonconstant stationary solution appears 
through a secondary bifurcation if $\Omega$ is symmetric (see also \cite{HV84}).
In~\cite{Mo90},
the stability of nonconstant stationary solutions was investigated
by constructing an invariant manifold (see also \cite{F90,MEF91,MJ92}).

\begin{figure}[htbp]
\centering
\includegraphics[keepaspectratio, scale=1.0]{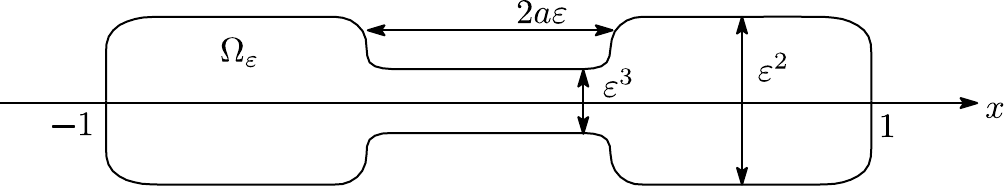}
\caption{Thin dumbbell-shaped domain.}
\label{tdsd}
\end{figure}

The problem \eqref{leq} is derived as another type of limiting equation.
For small $\vep>0$,
let $\Om_\vep$ be a thin dumbbell-shaped domain shown in Figure~\ref{tdsd}.
Then we can expect that solutions of \eqref{rde}
are approximated by those of some equation in one space dimension,
since $\Om_\vep$ shrinks to the interval $(-1,1)$ as $\vep \to 0$.
It is indeed shown in \cite{K} that stationary solutions of \eqref{rde} are approximated by 
solutions of \eqref{leq} in the following sense:
for any nondegenerate solution $u$ of \eqref{leq},
there exists a corresponding stationary solution $u_\vep$ of \eqref{rde} with $\Om=\Om_\vep$
such that $u_\vep$ converges to $u$ in an appropriate sense as $\vep \to 0$.
Therefore the analysis of \eqref{leq} provides 
information on the structure of solutions of \eqref{rde}.

Not only equations in higher dimensions,
but also other equations relate to \eqref{leq}.
Let $u$ be a solution of \eqref{leq} and let $\overline{u}$ be defined by
\begin{equation*}
\overline{u}(x):=\left\{
\begin{aligned}
&u(x+a) && (x \in [-1-a,-a)), 
\\
&\frac{u_x(+0)+u_x(-0)}{2}x +\frac{u(+0)+u(-0)}{2} && (x \in [-a,a]), 
\\
&u(x-a) && (x \in (a,1+a]).
\end{aligned}
\right.
\end{equation*}
Then, in a suitable sense, $\overline{u}$ satisfies the problem
\begin{equation*}
\left\{
\begin{aligned}
&(\overline{u})_{xx}+\lam (1-\chi_{(-a,a)}(x)) f(\overline{u})=0, \quad x \in (-1-a,1+a), 
\\
&(\overline{u})_x(-1-a)=(\overline{u})_x(1+a)=0,
\end{aligned}
\right.
\end{equation*}
where $\chi_A$ denotes the indicator function of a set $A$.
Therefore \eqref{leq} can be regarded as one of boundary value problems
for equations with variable coefficients.
In many such problems, secondary bifurcations are observed.
For instance, the Neumann problem for the equation $(c(x)u_x)_x+\lam (u-u^3)=0$ 
with a piecewise constant function $c(x)$ was studied by \cite{HR85},
and the Dirichlet problem for the equation $u_{xx}+d(\lam,x,u)=0$ 
with $d(\lam,x,u)=|x|^\lam u^p$, $|x|^\lam e^u$, or $(1-\chi_{(-\lam,\lam)}(x))u^p$
was investigated by \cite{KST18,ST19,T13,T17}.

We remark that the last two conditions of \eqref{leq} also appear 
when we consider the Schr\"{o}dinger operators with $\del'$-interactions
(for recent studies see \cite{AN13,AN14,AG18,AG20,GO20}).

The aim of this paper is to find bifurcation points 
on branches of solutions bifurcating from $u=0$,
and to determine the Morse index of solutions on the branches.
We first observe that each branch consists of odd solutions or even solutions.
Then, under additional assumptions on $f$, 
we show that a secondary bifurcation occurs only on the branch of odd solutions.
Finally, we prove that the Morse index of an odd solution with $m$ zeros
changes from $m+1$ to $m$ through a secondary bifurcation.
In particular, we conclude that 
the Morse index of a monotone odd solution becomes zero after the secondary bifurcation.
This fact is consistent with the result in \cite{V83}.


\subsection{Main result}

To state our main result,
we set up notation to be used throughout the paper.
Set
\begin{equation*}
X_0:=\{ u \in C^2([-1,1] \setminus \{0\});
 u_{xx}|_{[-1,0)} \mbox{ and } u_{xx}|_{(0,1]} \mbox{ are uniformly continuous}\}.
\end{equation*}
If $u \in X_0$, then $u|_{[-1,0)}$ (resp. $u|_{(0,1]}$) 
can be extended to a $C^2$ function on $[-1,0]$ (resp. $[0,1]$).
Hence we see that limits $u(\pm 0)$ and $u_x(\pm 0)$ exist for any $u \in X_0$
and that $X_0$ is a Banach space endowed with the norm
\begin{equation*}
\| u\|_{X_0}=\sup_{x \in [-1,1] \setminus \{0\}} |u(x)|
 +\sum_{k=1}^2 \sup_{x \in [-1,1] \setminus \{0\}} \left| \frac{d^k u}{dx^k}(x)\right|.
\end{equation*}
We focus on solutions of \eqref{leq} in the set $X \subset X_0$ defined by
\begin{equation*}
X:=\{ u \in X_0; |u(x)|<1 \mbox{ for } x \in [-1,1] \setminus \{0\} \}.
\end{equation*}
Let ${\cS}$ denote the set of all pairs 
$(\lam,u) \in (0,\infty) \times X$ satisfying \eqref{leq}:
\begin{equation*}
{\cS}:=\bigcup_{\lam \in (0,\infty)} \{ \lam\} \times {\cS}_\lam,
\qquad 
{\cS}_\lam :=\{ u \in X; u \mbox{ satisfies } \eqref{leq}\}.
\end{equation*}
Note that ${\cS}$ contains a trivial branch ${\cL}:=\{ (\lam,0)\}_{\lam \in (0,\infty)}$.
For ${\cA} \subset (0,\infty) \times X$,
we write
\begin{equation*}
-{\cA}:=\{ (\lam,-u); (\lam,u) \in {\cA} \}.
\end{equation*}
Then we see from the last condition of \eqref{basf} 
that $-{\cA} \subset {\cS}$ if ${\cA} \subset {\cS}$.
We define ${\cS}^o \subset {\cS}$ and ${\cS}^e \subset {\cS}$ by
\begin{gather*}
{\cS}^o:=\bigcup_{\lam \in (0,\infty)} \{ \lam\} \times {\cS}_\lam^o,
\qquad
{\cS}^e:=\bigcup_{\lam \in (0,\infty)} \{ \lam\} \times {\cS}_\lam^e,
\\
{\cS}_\lam^o:=\{ u \in {\cS}_\lam \setminus \{0\}; 
 u(-x)=-u(x) \mbox{ for } x \in [-1,1] \setminus \{0\} \},
\\
{\cS}_\lam^e:=\{ u \in {\cS}_\lam \setminus \{0\}; 
 u(-x)=u(x) \mbox{ for } x \in [-1,1] \setminus \{0\} \}.
\end{gather*}

For a solution $u \in {\cS}_\lam$, 
we also discuss the linearized eigenvalue problem
\begin{equation}\label{llevp}
\left\{
\begin{aligned}
&\varphi_{xx}+\lam f'(u)\varphi=\mu \varphi, \quad x \in (-1,1) \setminus \{ 0 \}, 
\\
&\varphi_x(-1)=\varphi_x(1)=0, 
\\
&\varphi(-0)+a\varphi_x(-0)=\varphi(+0)-a\varphi_x(+0),
\\
&\varphi_x(-0)=\varphi_x(+0).
\end{aligned}
\right.
\end{equation}
It is shown that the set of eigenvalues of \eqref{llevp} 
consists of a sequence of real numbers which diverges to $-\infty$
(see Lemma~\ref{lem:ipmi} in Section~\ref{sec:pre}). 
We say that $u$ is nondegenerate if all the eigenvalues are nonzero.
The number of positive eigenvalues is called the Morse index 
and is denoted by $i(u)$.


For $n \in \bbn$,
we define $\lam_n \in (0,\infty)$ by
\begin{equation*}
\lam_n:=\left\{
\begin{aligned}
&\frac{z_{k}^2}{f'(0)} && \mbox{if } n=2k-1, k \in \bbn,
\\
&\frac{(k\pi)^2}{f'(0)} && \mbox{if } n=2k, k \in \bbn,
\end{aligned}
\right.
\end{equation*}
where $z_k$ denotes the unique root of the equation $az\tan z=1$ in $((k-1)\pi,(k-1/2)\pi)$.
Note that $\lam_n<\lam_{n+1}$.
We set
\begin{equation*}
F(u):=2\int_0^u f(s)ds,
\end{equation*}
and fix $u_0 \in (0,1)$.

The main result of this paper is stated as follows.
\begin{thm}\label{mthm}
In addition to \eqref{basf}, assume that the following conditions are satisfied:
\begin{subequations}\label{aasfs}
\begin{gather}
f'(u)\left\{
\begin{aligned}
&>0 && \mbox{if } u \in (-u_0,u_0),
\\
&<0 && \mbox{if } u \in (-1,-u_0) \cup (u_0,1),
\end{aligned}
\right.
\label{aasfs0}
\\
\sgn (u) \frac{f'(u)F(u)}{f(u)^2} \mbox{ is strictly decreasing in }
 (-u_0,u_0) \setminus \{0\},
\label{aasfs1}
\\
\sgn (u) \frac{d}{du} \left( \frac{f'(u)F(u)^{\frac{3}{2}}}{f(u)^3}\right) \le 0
 \quad \mbox{for } u \in (-1,-u_0) \cup (u_0,1).
\label{aasfs2}
\end{gather}
\end{subequations}
Then for every $n \in \bbn$, 
there exists ${\cC}_n \subset {\cS}$
with the following properties.
\begin{enumerate}[(i)]
\item
${\cC}_n$ is a $C^1$ curve in $(0,\infty) \times X_0$ 
parametrized by $\lam \in (\lam_n,\infty)$ and bifurcates from $(\lam_n,0)$.
Moreover,
\begin{equation*}
{\cS}^o=\bigcup_{k=1}^\infty {\cC}_{2k-1} \cup (-{\cC}_{2k-1}), 
\qquad 
{\cS}^e=\bigcup_{k=1}^\infty {\cC}_{2k} \cup (-{\cC}_{2k}).
\end{equation*}
\item
There is no bifurcation point on ${\cL} \setminus \{ (\lam_n,0)\}_{n \in \bbn}$ and ${\cS}^e$.
\item
For every odd number $n \in \bbn$,
the curve ${\cC}_n$ (resp. $-{\cC}_n$) has a unique bifurcation point 
$(\lam^*_n,u^*_n)$ (resp. $(\lam^*_n,-u^*_n)$).
Furthermore, bifurcating solutions form a $C^1$ curve 
in a neighborhood of each bifurcation point.
\item
Let $(\lam,u) \in {\cC}_n \cup (-{\cC}_n)$.
If $n$ is even, then $u$ is nondegenerate and $i(u)=n$;
if $n$ is odd, then $u$ is nondegenerate unless $\lam=\lam^*_n$ and
\begin{equation*}
i(u)=\left\{
\begin{aligned}
&n && (\lam<\lam^*_n),
\\
&n-1 && (\lam \ge \lam^*_n).
\end{aligned}
\right.
\end{equation*}
Here $\lam^*_n$ is the number given in (iii).
\end{enumerate}
\end{thm}

\begin{figure}[h]
\includegraphics[scale=1.0]{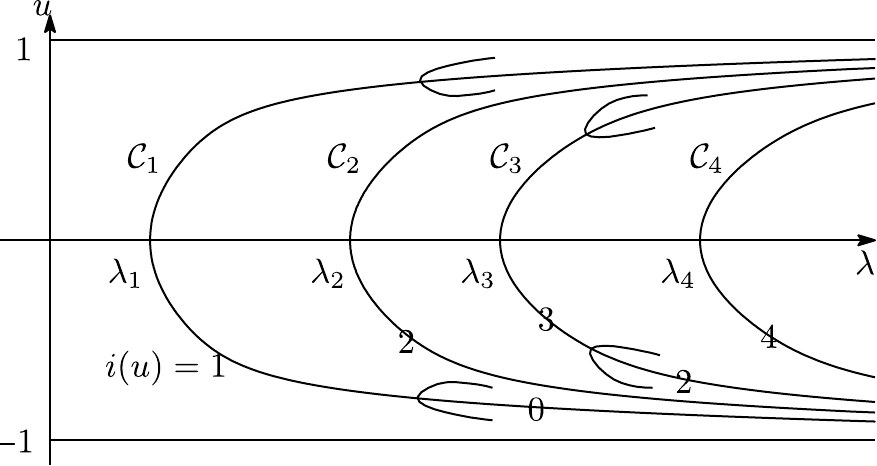}
\caption{Bifurcation diagram of \eqref{leq}.}
\label{bcdg}
\end{figure}

The bifurcation diagram of \eqref{leq} is shown in Figure~\ref{bcdg}.

\begin{remk}
$f(u)=u-u^3$ and $f(u)=\sin \pi u$ 
are typical examples satisfying \eqref{aasfs}, 
namely all of the conditions \eqref{aasfs0}, \eqref{aasfs1} and \eqref{aasfs2}.

The assumption \eqref{aasfs} is used 
to show the uniqueness of bifurcation points on ${\cC}_{2k-1}$.
The assertions (i) and (ii) are verified under the weaker assumption 
\begin{equation}\label{aasfw} 
\frac{f'(u)u}{f(u)}<1 \quad \mbox{for } u \in (-1,1) \setminus \{0\}.
\end{equation}
The fact that \eqref{aasfs} implies \eqref{aasfw} is not difficult to check.
For the convenience of the reader, 
we give the proof of this fact in Appendix~\ref{appendixA}.

Let $n$ be even and let $(\lam,u) \in {\cC}_n$. 
Then we have $u_x(+0)=u_x(-0)=0$, since $u$ is even.
Combining this with \eqref{leq} shows that $u$ can be extended smoothly up to $x=0$.
Therefore $u$ coincides with a solution $\tilde u$ 
of the usual Neumann problem for the equation
$\tilde u_{xx}+\lam f(\tilde u)=0$ in $(-1,1)$.
It is, however, not obvious that $u$ and $\tilde u$ have the same Morse index,
since the corresponding linearized eigenvalue problems are different.
\end{remk}

The main task in the proof of Theorem~\ref{mthm} 
is to investigate the behavior of eigenvalues of \eqref{llevp}.
For $(\lam,u) \in {\cC}_{2k-1}$ and $n \in \bbn \cup \{0\}$,
let $\mu_n (\lam)$ denote the $(n+1)$-th largest eigenvalue of \eqref{llevp}.
Then the condition $\mu_n (\lam_0)=0$ is necessary 
in order for a bifurcation to occur at $(\lam_0,u)$.
A theory of local bifurcations shows that $(\lam_0,u)$ is indeed a bifurcation point
if the additional condition $d\mu_n(\lam_0)/d\lam \neq 0$ is satisfied.
The main difficulty arises from the verification of this condition.
We will show in Sections~\ref{sec:sb} and \ref{sec:pt} that the assumption \eqref{aasfs} 
gives $d\mu_{2k-2}(\lam^*)/d\lam<0$ for any $\lam^*$ satisfying $\mu_{2k-2} (\lam^*)=0$.
This enables us to conclude that a secondary bifurcation occurs on ${\cC}_{2k-1}$ just once.

This paper is organized as follows.
Section~\ref{sec:pre} provides some preliminaries.
A main task is to reduce the problem \eqref{leq} to a finite dimensional problem.
In Section~\ref{sec:pb}, 
we study primary branches bifurcating from the trivial branch $\cL$.
In Section~\ref{sec:sb},
we discuss the number of bifurcation points on the primary branches.
Section~\ref{sec:pt} is devoted to the evaluation of the Morse index
and the proof of Theorem~\ref{mthm}.

\section{Preliminaries}\label{sec:pre}

In this section,
we first introduce a function $G$ to be used throughout the study,
next convert \eqref{leq} into a finite dimensional problem by the shooting method,
then study general properties of eigenpairs of \eqref{llevp},
and finally give sufficient conditions for a bifurcation 
by applying a bifurcation theorem developed in \cite{CR71}.

\subsection{Definition and properties of $G$}

Set 
\begin{equation*}
\be_0:=\sqrt{2\int_0^1 f(s)ds} \left( =\sqrt{F(1)}=\sqrt{F(-1)} \right),
\qquad
I:=(-\be,\be_0).
\end{equation*}
From \eqref{basf}, we see that the function
\begin{equation}\label{Ginv}
(-1,1) \ni u \mapsto \sgn (u) \sqrt{F(u)} \in I
\end{equation}
is strictly increasing.
We then define a function $G:I \to (-1,1)$ to be the inverse of this function,
that is, $G(v)$ is determined by the relation
\begin{equation*}
F(G(v))=v^2, \quad |G(v)|<1, \quad G(v) \left\{ 
\begin{aligned}
&<0 && \mbox{if } v \in (-\be_0,0), \\
&>0 && \mbox{if } v \in (0,\be_0).
\end{aligned}
\right.
\end{equation*}
Furthermore, we set 
\begin{equation}\label{hdef}
h(v):= 
1-\frac{G''(v)v}{G'(v)},
\qquad
H(v):=v^2 -\frac{G(v)v}{G'(v)}.
\end{equation}

In the following lemma,
we collect properties and formulas for $G$.
Although part of the lemma is shown in \cite{S90},
we prove all the assertions in Appendix~\ref{appendixB} for readers' convenience.
\begin{lem}\label{lem:Gpro}
The following hold.
\begin{enumerate}[(i)]
\item
$G(v) \in C^2(I) \cap C^3(I \setminus \{0\})$
and $G''(v)v \in C^1(I)$.
\item
\eqref{aasfs0}, \eqref{aasfs1} and \eqref{aasfs2} 
are equivalent to the following conditions 
\eqref{aashs0}, \eqref{aashs1} and \eqref{aashs2}, respectively:
\begin{subequations}\label{aashs}
\begin{gather}
h(v)\left\{
\begin{aligned}
&>0 && \mbox{if } v \in (-v_0,v_0),
\\
&<0 && \mbox{if } v \in (-\be_0,-v_0) \cup (v_0,\be_0), 
\end{aligned}
\right.
\label{aashs0} 
\\
\sgn (v) h(v) \mbox{ is strictly decreasing in } (-v_0,v_0) \setminus \{0\},
\label{aashs1} 
\\
\sgn (v) \frac{d}{dv} (G'(v)h(v)) \le 0 \quad \mbox{for } v \in (-\be_0,-v_0) \cup (v_0,\be_0).
\label{aashs2} 
\end{gather}
\end{subequations}
Here $v_0=G^{-1}(u_0)$. 
Moreover, \eqref{aasfw} holds if and only if
\begin{equation}\label{aasGw}
\sgn(v) H'(v)>0 \quad \mbox{for } v \in I \setminus \{0\}.
\end{equation}
\item
There are constants $c>0$ and $C>0$ such that for all $v \in I$,
\begin{gather}
\frac{c}{\sqrt{\be_0-|v|}} \le G'(v) \le \frac{C}{\sqrt{\be_0-|v|}}, 
\label{Gpas} \\
\sgn(v) G''(v) \ge \frac{c}{(\be_0-|v|)^{3/2}}-C.
\label{Gdpas} 
\end{gather}
\end{enumerate}
\end{lem}

The function $G$ is introduced in \cite{S90,W83,W86}
to obtain a simple expression of solutions of the equation $w_{xx}+f(w)=0$.
A solution $w$ satisfying $|w|<1$
corresponds to the closed orbit $w_x^2+F(w)=c_0$ in the $ww_x$-plane,
where $c_0$ is a nonnegative constant. 
By the change of variables $w=G(\tilde w)$,
the orbit is transformed into the circle $w_x^2+\tilde w^2=c_0$ in the $\tilde ww_x$-plane.
Hence $w$ is written as $w(x)=G(\sqrt{c_0} \cos \rho (x))$
for some suitable function $\rho (x)$.
The details of this argument are given in the next subsection.

\subsection{Reduction to a finite dimensional problem}

For $\be_1,\be_2 \in I$ and $\lam \in (0,\infty)$,
let $u_1$ and $u_2$ to be solutions of the initial value problems 
\begin{equation}\label{u12ivp}
\left\{
\begin{aligned}
&(u_1)_{xx}+\lam f(u_1)=0, \quad x \in \bbr, \\
&u_1(-1)=G(\be_1), \ (u_1)_x(-1)=0,
\end{aligned}
\right.
\qquad \qquad
\left\{
\begin{aligned}
&(u_2)_{xx}+\lam f(u_2)=0, \quad x \in \bbr, \\
&u_2(1)=G(\be_2), \ (u_2)_x(1)=0.
\end{aligned}
\right.
\end{equation}
Then the function $u$ defined by
\begin{equation*}
u(x)=\left\{
\begin{aligned}
&u_1(x) && (x \in [-1,0)), \\
&u_2(x) && (x \in (0,1])
\end{aligned}
\right.
\end{equation*}
satisfies \eqref{leq} if and only if
\begin{equation}\label{smoeq}
\left\{
\begin{aligned}
u_1(0)+a(u_1)_x(0)&=u_2(0)-a(u_2)_x(0), \\
(u_1)_x(0)&=(u_2)_x(0).
\end{aligned}
\right.
\end{equation}
Now we introduce a solution $(U,V)=(U(y,\be),V(y,\be))$ of the initial value problem
\begin{equation}\label{UVeq}
\left\{
\begin{aligned}
&U_y=V, \quad V_y=-f(U), && y \in \bbr, \\
&(U(0,\be),V(0,\be))=(G(\be),0).
\end{aligned}
\right.
\end{equation}
Then we have
\begin{gather}
u_1(x)=U\left( \sqrt{\lambda} (x+1),\be_1 \right),
\quad
u_2(x)=U\left( \sqrt{\lambda} (1-x),\be_2 \right),
\label{u12U}
\\
(u_1)_x (x)=\sqrt{\lambda} V\left( \sqrt{\lambda} (x+1),\be_1 \right),
\quad
(u_2)_x (x)=-\sqrt{\lambda} V\left( \sqrt{\lambda} (1-x),\be_2 \right).
\nonumber
\end{gather}
Hence \eqref{smoeq} is equivalent to the equation
\begin{equation}\label{smteq}
(P(\lam,\be_1),Q(\lam,\be_1))=(P(\lam,\be_2),-Q(\lam,\be_2)),
\end{equation}
where
\begin{equation*}
P(\lam,\be):=U(\sqrt{\lambda},\be) +a\sqrt{\lambda} V(\sqrt{\lambda},\be), 
\qquad
Q(\lam,\be):=V(\sqrt{\lambda},\be).
\end{equation*}
We remark that 
\begin{gather}
U(x,-\be)=-U(x,\be),
\qquad
V(x,-\be)=-V(x,\be),
\label{UVsy}
\\
P(\lam,-\be)=-P(\lam,\be),
\qquad
Q(\lam,-\be)=-Q(\lam,\be),
\label{PQsy}
\end{gather}
which follow from the fact that $f$ and $G$ are odd functions.

Let us investigate the relation between \eqref{leq} and \eqref{smteq}.
The set of solutions of \eqref{smteq} is denoted by 
\begin{gather*}
{\cT}:=\bigcup_{\lam \in (0,\infty)} \{ \lam\} \times {\cT}_\lam,
\\
\mathcal{T}_\lam :=\{ (\be_1,\be_2) \in I \times I;
 (P(\lam,\be_1),Q(\lam,\be_1))=(P(\lam,\be_2),-Q(\lam,\be_2)) \}.
\end{gather*}
We define $C^2$ mappings $\Phi :X \to I \times I$ and $\Psi_\lam :I \times I \to X$ by
\begin{gather*}
\Phi (u):=(G^{-1}(u(-1)),G^{-1}(u(1))),
\\
\Psi_\lam (\be_1,\be_2)(x) :=\left\{
\begin{aligned}
&U\left( \sqrt{\lambda} (x+1),\be_1 \right) (=u_1(x)) && \mbox{for } x \in [-1,0), \\
&U\left( \sqrt{\lambda} (1-x),\be_2 \right) (=u_2(x)) && \mbox{for } x \in (0,1].
\end{aligned}
\right.
\end{gather*}
We note that
\begin{equation}\label{Pssy}
\Psi_\lam (-\be_1,-\be_1)(x)=-\Psi_\lam (\be_1,\be_2)(x),
\end{equation}
which is obtained immediately by \eqref{UVsy}.

By the following lemma, we see that
\eqref{leq} is reduced to the finite dimensional problem \eqref{smteq}.

\begin{lem}\label{ppoo}
$\Phi|_{\mathcal{S}_\lam}$ and $\Psi_\lam|_{\mathcal{T}_\lam}$ 
are one-to-one correspondences between $\mathcal{S}_\lam$ and $\mathcal{T}_\lam$.
More precisely, 
$\Phi(\mathcal{S}_\lam)=\mathcal{T_\lam}$, $\Psi_\lam(\mathcal{T}_\lam)=\mathcal{S}_\lam$,
$\Psi_\lam \circ \Phi|_{\mathcal{S}_\lam}=id_{\mathcal{S}_\lam}$ 
and $\Phi \circ \Psi_\lam|_{\mathcal{T}_\lam}=id_{\mathcal{T}_\lam}$.
\end{lem}

\begin{proof}
That $\Psi_\lam(\mathcal{T}_\lam) \subset \mathcal{S}_\lam$ 
and $\Phi \circ \Psi_\lam|_{\mathcal{T}_\lam}=id_{\mathcal{T}_\lam}$ 
follow from the definitions of $\Phi$ and $\Psi_\lam$.
To see $\Phi(\mathcal{S}_\lam) \subset \mathcal{T}_\lam$ and
$\Psi_\lam \circ \Phi|_{\mathcal{S}_\lam}=id_{\mathcal{S}_\lam}$,
let $u \in \mathcal{S}_\lam$ and $(\be_1,\be_2)=\Phi(u)$.
Then, by the first two conditions of \eqref{leq}, 
we see that $u_1:=u|_{[-1,0)}$ and $u_2:=u|_{(0,1]}$ satisfy \eqref{u12ivp}.
Hence $u_1$ and $u_2$ are given by \eqref{u12U},
which gives $(\Psi_\lam \circ \Phi)(u)=u$.
Since the last two conditions of \eqref{leq} lead to \eqref{smoeq},
we deduce that $\Phi (u)=(\be_1,\be_2) \in \mathcal{T}_\lam$.
We have thus proved $\Phi(\mathcal{S}_\lam) \subset \mathcal{T}_\lam$ and
$\Psi_\lam \circ \Phi|_{\mathcal{S}_\lam}=id_{\mathcal{S}_\lam}$.
This completes the proof.
\end{proof}

The nondegeneracy of a solution of \eqref{leq} 
can also be determined by the corresponding solution of \eqref{smteq}.
We set
\begin{align*}
D(\lam,\be_1,\be_2):=&\det 
\begin{pmatrix}
P_\be (\lam,\be_1) & -P_\be (\lam,\be_2) \\
Q_\be (\lam,\be_1) & Q_\be (\lam,\be_2)
\end{pmatrix}
\\
=&P_\be (\lam,\be_1)Q_\be (\lam,\be_2) +Q_\be (\lam,\be_1)P_\be (\lam,\be_2),
\end{align*}
where $P_\be$ and $Q_\be$ stand for derivatives with respect to $\be$.
\begin{lem}\label{lem:ndcD}
Let $u \in \mathcal{S}_\lam$ and put $(\be_1,\be_2)=\Phi (u) \in {\mathcal T}_\lam$.
Then $u$ is nondegenerate if and only if $D(\lam,\be_1,\be_2) \neq 0$.
\end{lem}

\begin{proof}
Let $u_1$ and $u_2$ be given by \eqref{u12U} and define
\begin{equation}\label{phjdef}
\phi_j:=\frac{1}{G'(\be_j)} \frac{\D u_j}{\D \be_j},
\quad
j=1,2.
\end{equation}
Then, by the definitions of $P$ and $Q$, we have 
\begin{equation}\label{PQphj}
P_\be (\lam,\be_j)=G'(\be_j) \left. \left\{ \phi_j +(-1)^{j-1} a(\phi_j)_x\right\} \right|_{x=0},
\quad
Q_\be (\lam,\be_j)=\frac{(-1)^{j-1}G'(\be_j)}{\sqrt{\lam}} (\phi_j)_x |_{x=0}.
\end{equation}
Hence the condition $D(\lam,\be_1,\be_2)=0$ holds if and only if 
\begin{equation}\label{Ddcon}
\mbox{there is } (\al_1,\al_2) \neq (0,0) \mbox{ such that } 
\left\{
\begin{aligned}
\al_1 \{ \phi_1+a(\phi_1)_x\}|_{x=0}&=\al_2 \{ \phi_2-a(\phi_2)_x\}|_{x=0}, &&
\\
\al_1 (\phi_1)_x|_{x=0}&=\al_2 (\phi_2)_x|_{x=0},
&&
\end{aligned}
\right.
\end{equation}
which means that $(P_\be(\lam,\be_1), Q_\be(\lam,\be_1))$ 
and $(-P_\be(\lam,\be_2), Q_\be(\lam,\be_2))$ are linearly dependent.

Suppose that $D(\lam,\be_1,\be_2)=0$, 
that is, the condition \eqref{Ddcon} holds.
By definition, we see that $\phi_j$ satisfies 
\begin{equation}\label{duaeq}
(\phi_j)_{xx} +\lam f'(u_j) \phi_j=0, 
\qquad 
\phi_j|_{x=(-1)^j}=1,
\qquad
(\phi_j)_x|_{x=(-1)^j}=0.
\end{equation}
We define $\varphi \in X_0 \setminus \{0\}$ by
\begin{equation*}
\varphi (x)=\left\{
\begin{aligned}
&\al_1 \phi_1(x) && (x \in [-1,0)), \\
&\al_2 \phi_2(x) && (x \in (0,1]).
\end{aligned}
\right.
\end{equation*}
Then we see from \eqref{Ddcon} and \eqref{duaeq} that 
$\varphi$ satisfies \eqref{llevp} for $\mu=0$.
This means that $u$ is not nondegenerate.

Suppose conversely that $u$ is not nondegenerate, that is,
there is $\varphi \in X_0 \setminus \{0\}$ satisfying \eqref{llevp} for $\mu=0$.
If $\varphi$ vanished at $x=-1$ or $x=1$, 
we would have $\varphi \equiv 0$
due to the uniqueness of solutions of initial value problems. 
Hence both $\varphi (-1)$ and $\varphi (1)$ are nonzero.
Since $\phi_j$ and $\varphi$ satisfy the same differential equation,
we deduce that
\begin{equation}\label{pjvpj}
\phi_1 (x)=\frac{1}{\varphi (-1)}\varphi (x) \quad (x \in [-1,0)), 
\quad 
\phi_2 (x)=\frac{1}{\varphi (1)}\varphi (x) \quad (x \in (0,1]).
\end{equation}
This implies that \eqref{Ddcon} holds 
for $(\al_1,\al_2)=(\varphi (-1),\varphi (1)) \neq (0,0)$.
Thus $D(\lam,\be_1,\be_2)=0$, and the lemma follows.
\end{proof}

Let us derive explicit formulas for $P$ and $Q$ by solving \eqref{UVeq}.
For $(y,\be) \in \bbr \times I$, 
let $\Theta \in \bbr$ be given implicitly by the relation
\begin{equation}\label{Thdef}
\int_0^\Theta G'(\be \cos \tau)d\tau =y.
\end{equation}
Since $G'$ is positive, 
the left-hand side is strictly increasing with respect to $\Theta$
and diverges to $\pm \infty$ as $\Theta \to \pm \infty$.
Hence $\Theta=\Theta (y,\be)$ 
is uniquely determined for every $(y,\be) \in \bbr \times I$.
Moreover, by the implicit function theorem,
we infer that $\Theta$ is of class $C^1$.
We put $\theta (\lam,\be) :=\Theta (\sqrt{\lam},\be)$,
that is, $\theta$ is determined by
\begin{equation}\label{tblr}
\int_0^\theta G'(\be \cos \tau)d\tau =\sqrt{\lambda}.
\end{equation}
We note that $\theta>0$.

\begin{lem}
$P$ and $Q$ are written as
\begin{equation}\label{PQrep}
P=G(\be \cos \theta)-a\sqrt{\lambda} \be \sin \theta, 
\qquad
Q=-\be \sin \theta.
\end{equation}
Furthermore, 
\begin{equation}\label{pqdb}
\left\{
\begin{aligned}
P_\be&=G'(\be \cos \theta) \cos \theta
 -a\sin \theta \int_0^\theta G'(\be \cos \tau) d\tau
\\
&\quad +\left( \be \sin \theta +a\frac{\be \cos \theta}{G'(\be \cos \theta)}
 \int_0^\theta G'(\be \cos \tau) d\tau \right)
  \int_0^\theta G''(\be \cos \tau) \cos \tau d\tau,
\\
Q_\be&=-\sin \theta+\frac{\be \cos \theta}{G'(\be \cos \theta)}
 \int_0^\theta G''(\be \cos \tau) \cos \tau d\tau.
\end{aligned}
\right.
\end{equation}
\end{lem}

\begin{proof}
Define
\begin{equation*}
(U,V):=(G(\be \cos \Theta),-\be \sin \Theta).
\end{equation*}
Differentiating both sides of \eqref{Thdef} 
yields $G'(\be \cos \Theta) \Theta_y=1$,
and hence
\begin{gather*}
U_y=G'(\be \cos \Theta) \cdot (-\be \sin \Theta) \cdot \Theta_y=V,
\\
V_y=(-\be \cos \Theta) \cdot \Theta_y
 =-\frac{\be \cos \Theta}{G'(\be \cos \Theta)}=-f(U),
\end{gather*}
where we have used \eqref{dGfo} in deriving the last equality.
Since $\Theta|_{y=0}=0$, we have $(U,V)|_{y=0}=(G(\be),0)$.
This shows that $(U,V)$ satisfies \eqref{UVeq},
and therefore we obtain \eqref{PQrep}.

By \eqref{PQrep},
we have
\begin{equation}\label{pqdb0}
\left\{
\begin{aligned}
&P_\be=G'(\be \cos \theta)(\cos \theta -\be \theta_\be \sin \theta)
 -a\sqrt{\lam} (\sin \theta +\be \theta_\be \cos \theta),
\\
&Q_\be=-\sin \theta -\be \theta_\be \cos \theta.
\end{aligned}
\right.
\end{equation}
Differentiating \eqref{tblr} gives
\begin{equation}\label{thdb}
\theta_\be =-\frac{1}{G'(\be \cos \theta))} \int_0^\theta G''(\be \cos \tau) \cos \tau d\tau.
\end{equation}
\eqref{pqdb} is then derived easily
by plugging \eqref{tblr} and \eqref{thdb} into \eqref{pqdb0}.
\end{proof}

\subsection{Properties of eigenpairs of \eqref{llevp}}

We recall known facts for the eigenvalue problem
\begin{equation}\label{gevp}
\left\{
\begin{aligned}
&\psi_{xx}+q(x) \psi=\nu w(x) \psi, && x \in (b,c), \\
&\psi_x(b)=\psi_x(c)=0,
\end{aligned}
\right.
\end{equation}
where $q$ and $w$ are given functions and $b,c \in \bbr$, $b<c$.
Since we deal with the case where $q$ and $w$ are not necessarily continuous,
we consider eigenfunctions in a generalized sense:
by an eigenfunction of \eqref{gevp} 
we mean a function $\psi \in W^{2,1}(b,c) \setminus \{0\}$ 
satisfying the differential equation in \eqref{gevp} almost everywhere
and the boundary condition in the usual sense
(note that $W^{2,1}(b,c) \subset C^1([b,c])$).

\begin{thm}[\cite{A64}]\label{thm:gevp}
Suppose that $q$ and $w$ satisfy
\begin{equation}\label{assmp:gevp}
\left\{
\begin{aligned}
&q, w \in L^1(b,c),
\\
&w \ge 0 \quad \mbox{in } (b,c),
\\
&w>0 \quad \mbox{in } (b,b+\del) \cup (c-\del,c) \mbox{ for some } \del>0,
\\
&\{ w=0\} \subset \{ q=0\}.
\end{aligned}
\right.
\end{equation}
Then the following hold.
\begin{enumerate}[(i)]
\item
The eigenvalues of \eqref{gevp} 
consists of real numbers $\{ \nu_n\}_{n=0}^\infty$ satisfying
\begin{equation*}
\nu_n>\nu_{n+1},
\qquad
\nu_n \to -\infty
\quad (n \to \infty).
\end{equation*}
Furthermore, each eigenvalue is simple, that is,
the eigenspace associated with $\nu_n$ is one-dimensional.
\item
Any eigenfunction corresponding to $\nu_n$ have exactly $n$ zeros in $(b,c)$. 
\end{enumerate}
\end{thm}

For the proof of this theorem, see \cite[Theorems~8.4.5 and 8.4.6]{A64}.

In what follows, 
we denote by $\cE$ the set of all pairs $(q,w)$ satisfying \eqref{assmp:gevp},
and write $\nu_n(q)$ instead of $\nu_n$ if we emphasize the dependence on $q$.

\begin{lem}\label{lem:coev}
For $n \in \bbn \cup \{0\}$, the following hold.
\begin{enumerate}[(i)]
\item
Let $(q,w), (\tilde q,w) \in {\cE}$.
Suppose that $q \ge \tilde q$ in $(b,c)$ 
and $q>\tilde q$ in some nonempty open subinterval of $(b,c)$. 
Then $\nu_n(q)>\nu_n(\tilde q)$.
\item
For any $\vep>0$ and $(q,w) \in {\cE}$, 
there exists $\del>0$ such that $|\nu_n(q)-\nu_n(\tilde q)|<\vep$ 
whenever $(\tilde q,w) \in {\cE}$ and $\| q -\tilde q\|_{L^1(b,c)}<\del$.
In other words, the mapping $q \mapsto \nu_n (q)$ is continuous.
\end{enumerate}
\end{lem}

In the case where $w$ is positive,
the above lemma is well-known 
and can be shown by an argument based on Pr\"{u}fer's transformation.
It is not difficult to check that the same argument works 
under the conditions stated in the lemma.
We omit the detailed proof.

Let us examine the properties of eigenvalues of \eqref{llevp}
using Theorem~\ref{thm:gevp} and Lemma~\ref{lem:coev}.
For $u \in X_0$,
we define a function $\overline{u}$ by
\begin{equation*}
\overline{u}(x):=\left\{
\begin{aligned}
&u(x+a) && (x \in [-1-a,-a)), 
\\
&\frac{u_x(+0)+u_x(-0)}{2}x +\frac{u(+0)+u(-0)}{2} && (x \in [-a,a]), 
\\
&u(x-a) && (x \in (a,1+a]).
\end{aligned}
\right.
\end{equation*}
Then it is straightforward to show that 
if $u(-0)+au_x(-0)=u(+0)-au_x(+0)$ and $u_x(-0)=u_x(+0)$,
then $\overline{u} \in W^{2,\infty}(-1-a,1+a)$.
We set
\begin{equation*}
w_0(x):=\left\{ 
\begin{aligned}
&1 && (x \in [-1-a,-a) \cup (a,1+a]), 
\\
&0 && (x \in [-a,a]).
\end{aligned}
\right.
\end{equation*}

\begin{lem}\label{lem:ipmi}
Let $u \in X_0$.
Then the set of eigenvalues of \eqref{llevp} 
consists of an infinite sequence of real numbers which diverges to $-\infty$.
Furthermore, each eigenvalue is simple and continuous 
with respect to $(\lam,u) \in (0,\infty) \times X_0$.
\end{lem}

\begin{proof}
Assume that $(\mu,\varphi)$ is an eigenpair of \eqref{llevp}.
Then $\overline{\varphi}$ satisfies $\overline{\varphi} \in W^{2,\infty}(-1-a,1+a)$,
$(\overline{\varphi})_x(-1-a)=(\overline{\varphi})_x(1+a)=0$ and 
\begin{equation*}
(\overline{\varphi})_{xx} +\lam f'(\overline{u}) \overline{\varphi}=\mu \overline{\varphi}
 \quad \mbox{in } (-1-a,-a) \cup (a,1+a),
\qquad
(\overline{\varphi})_{xx}=0 \quad \mbox{in } (-a,a).
\end{equation*}
Hence we see that $(\mu,\overline{\varphi})$ is an eigenpair of \eqref{gevp} with
\begin{equation}\label{bcqw}
b=-1-a, \quad c=1+a, \quad q=\lam w_0 f'(\overline{u}), \quad w=w_0.
\end{equation}

Conversely, let $(\nu,\psi)$ be an eigenpair of \eqref{gevp}
for $b$, $c$, $q$, $w$ given above.
Then we have $\psi_{xx}=0$ in $(-a,a)$,
which shows that $\psi$ is a linear function in $(-a,a)$.
Therefore
\begin{equation}\label{pssr}
\psi_x (-a)=\psi_x (a)=\frac{\psi (a) -\psi (-a)}{2a}
 \big( =(\mbox{the slope of the graph of } \psi \mbox{ in } (-a,a))\big).
\end{equation}
We put
\begin{equation*}
\underline{\psi} (x)=\left\{
\begin{aligned}
&\psi (x-a) && (x \in [-1,0)), 
\\
&\psi (x+a) && (x \in (0,1]).
\end{aligned}
\right.
\end{equation*}
Since $q=\lam w_0 f'(\overline{u})$ and $w=w_0$ 
are uniformly continuous on $[-1-a,a) \cup (a,1+a]$,
we see that $\psi_{xx}(=-q\psi +\nu w\psi)$ is also uniformly continuous on the same region.
Hence we have $\underline{\psi} \in X_0$. 
Furthermore, \eqref{gevp} and \eqref{pssr} imply that 
\eqref{llevp} is satisfied for $(\mu,\varphi)=(\nu,\underline{\psi})$.
We thus conclude that $(\nu,\underline{\psi})$ is an eigenpair of \eqref{llevp}.

The above argument shows that the set of eigenvalues of \eqref{llevp}
coincides with that of \eqref{gevp}.
It is straightforward to check that 
$q=\lam w_0 f'(\overline{u})$ and $w=w_0$ satisfy \eqref{assmp:gevp},
and therefore we obtain the desired conclusion 
by (i) of Theorem~\ref{thm:gevp} and (ii) of Lemma~\ref{lem:coev}.
\end{proof}

From now on, for $n \in \bbn \cup \{0\}$,
let $\mu_n(u)$ stand for the $(n+1)$-th largest eigenvalue of \eqref{llevp}.
\begin{lem}\label{lem:veef}
Let $\varphi$ be an eigenfunction of \eqref{llevp} corresponding to $\mu_n(u)$.
Then 
\begin{equation*}
\varphi (-1)\varphi (1)
\left\{
\begin{aligned}
&>0 && \mbox{if } n \mbox{ is even},
\\
&<0 && \mbox{if } n \mbox{ is odd}.
\end{aligned}
\right.
\end{equation*}
\end{lem}

\begin{proof}
As shown in the proofs of Lemmas~\ref{lem:ndcD} and \ref{lem:ipmi},
we know that $\varphi (-1)\varphi (1)$ is nonzero
and that $(\mu_n(u),\overline{\varphi})$ is the $(n+1)$-th eigenpair of \eqref{gevp} 
with $b,c,q,w$ given by \eqref{bcqw}.
By (ii) of Theorem~\ref{thm:gevp},
we see that $\overline{\varphi}$ has exactly $n$ zeros in $(-1-a,1+a)$.
Hence $\overline{\varphi}(-1-a)=\varphi (-1)$ and $\overline{\varphi}(-1-a)=\varphi (1)$ 
have the same sign if $n$ is even and have opposite signs if $n$ is odd.
\end{proof}

We conclude this subsection with a lemma which provides basic estimates of the Morse index.
\begin{lem}\label{lem:Mies}
Suppose that \eqref{aasfw} holds and that $u \in {\cS}_\lam \setminus \{0\}$ 
vanishes at exactly $n$ points in $(-1,1) \setminus \{0\}$.
Then $\mu_{n+1}(u)<0$ if $u(-0)u(+0)<0$ and $\mu_{n}(u)<0$ if $u(-0)u(+0)>0$.
\end{lem}

\begin{proof}
Set
\begin{equation*}
q:=\left\{
\begin{aligned}
&\lam w_0 \frac{f(\overline{u})}{\overline{u}} && \mbox{if } \ \overline{u} \neq 0, 
\\
&\lam w_0 f'(0) && \mbox{if } \ \overline{u}=0,
\end{aligned}
\right.
\qquad
\tilde q:=\lam w_0 f'(\overline{u}).
\end{equation*}
Then one can easily check that $(q,w_0), (\tilde q,w_0) \in {\cE}$.
From the assumptions $u \not\equiv 0$ and \eqref{aasfw}, 
we see that $q \ge \tilde q$ in $(-1-a,1+a)$ 
and $q>\tilde q$ in a nonempty open subinterval of $(-1-a,1+a)$.
Hence (i) of Lemma~\ref{lem:coev} shows that $\nu_j (q)>\nu_j(\tilde q)$ for all $j$.
As shown in the proof Lemma~\ref{lem:ipmi},
we know that $\mu_j(u)=\nu_j(\tilde q)$.
Therefore the lemma is proved if we show that
\begin{equation}\label{tu0ev}
\nu_{n+1} (q)=0 \ \mbox{ if } \ u(-0)u(+0)<0,
\qquad
\nu_n (q)=0 \ \mbox{ if } \ u(-0)u(+0)>0.
\end{equation}

The assumption $u \in {\cS}_\lam \setminus \{0\}$ shows that
$\overline{u} \in W^{2,\infty}(-1-a,1+a) \setminus \{0\}$ and
\begin{equation*}
\left\{
\begin{aligned}
&(\overline{u})_{xx}+\lam w_0(x) f(\overline{u})=0, && x \in (-1-a,1+a), 
\\
&(\overline{u})_x(-1-a)=(\overline{u})_x(1+a)=0. 
\end{aligned}
\right.
\end{equation*}
Note that the above equation is written as 
$(\overline{u})_{xx}+q(x) \overline{u}=0$.
Hence we infer that $\nu_m(q)=0$ for some $m$
and $\overline{u}$ is an eigenfunction corresponding to $\nu_m(q)=0$.
Since $u$ is assumed to vanish at exactly $n$ points in $(-1,1) \setminus \{0\}$,
we deduce that
\begin{equation*}
(\mbox{the number of zeros of } \overline{u} \mbox{ in } (-1-a,1+a))
=\left\{
\begin{aligned}
&n+1 && \mbox{if } u(-0)u(+0)<0,
\\
&n && \mbox{if } u(-0)u(+0)>0.
\end{aligned}
\right.
\end{equation*}
By (ii) of Theorem~\ref{thm:gevp}, 
we conclude that $m=n+1$ if $u(-0)u(+0)<0$ and $m=n$ if $u(-0)u(+0)>0$.
Thus \eqref{tu0ev} is verified, and the proof is complete.
\end{proof}

\subsection{Conditions for a bifurcation}

In this subsection,
we observe that sufficient conditions for a solution to be a bifurcation point
are described by means of $D(\lam,\be_1,\be_2)$.
The precise statement is given as follows.
\begin{prop}\label{prop:bt}
Let $J \subset (0,\infty)$ be an open interval containing a point $\lam_0$
and suppose that ${\cC}=\{ (\lam,u(\cdot,\lam))\}_{\lam \in J} \subset {\cS}$ 
is a $C^1$ curve in $(0,\infty) \times X_0$.
Set $(\be_1(\lam),\be_2(\lam)):=\Phi (u(\cdot,\lam))$.
Then the following hold.
\begin{enumerate}[(i)]
\item
If $D(\lam_0,\be_1(\lam_0),\be_2(\lam_0)) \neq 0$, then 
there is a neighborhood 
${\cN}$ of $(\lam_0,u(\cdot,\lam_0))$ in $(0,\infty) \times X_0$
such that ${\cS} \cap {\cN}={\cC} \cap {\cN}$.
\item
Suppose that
\begin{equation}\label{assmp:bt}
D(\lam_0,\be_1(\lam_0),\be_2(\lam_0))=0,
\qquad
\left. \frac{d}{d\lam} D(\lam,\be_1(\lam),\be_2(\lam)) \right|_{\lam=\lam_0} \neq 0.
\end{equation}
Then there exists ${\tilde \cC} \subset {\cS}$
such that
\begin{gather*}
{\tilde \cC} \mbox{ is a } C^1 \mbox{ curve in } (0,\infty) \times X_0
\mbox{ intersecting } {\cC} \mbox{ transversally at } (\lam_0,u(\cdot,\lam_0)),
\\
{\cS} \cap {\cN}=({\cC} \cup {\tilde \cC}) \cap {\cN}
\mbox{ for some neighborhood } {\cN} \mbox{ of } (\lam_0,u(\cdot,\lam_0)) 
\mbox{ in } (0,\infty) \times X_0.
\end{gather*}
\item
Let \eqref{assmp:bt} hold, 
and take $n \in \bbn \cup \{0\}$ satisfying $\mu_n (u(\cdot,\lam_0))=0$.
Then $\mu_n (u(\cdot,\lam))$ is continuously differentiable 
in a neighborhood of $\lam_0$ and
\begin{equation}\label{evdf}
\sgn \left( \left. \frac{d}{d\lam} \mu_n (u(\cdot,\lam)) \right|_{\lam=\lam_0} \right)
 =(-1)^{n-1} \sgn \left( 
  \left. \frac{d}{d\lam} D(\lam,\be_1(\lam),\be_2(\lam)) \right|_{\lam=\lam_0} \right).
\end{equation}
\end{enumerate}
\end{prop}

\begin{remk}
The existence of the number $n$ appeared in (iii) is guaranteed by Lemma~\ref{lem:ndcD}.
\end{remk}

To show the above proposition, 
we recall an abstract bifurcation theorem proved in \cite[Theorem~1.7]{CR71}.
In what follows,
we denote by $K(T)$ and $R(T)$ the kernel and range of a linear operator $T$, respectively.
\begin{thm}[\cite{CR71}]\label{thm:crbt}
Let $\mathcal{X}$ and $\mathcal{Y}$ be Banach spaces, $J$ an open interval
and $\mathcal{F}=\mathcal{F}(\lam,w)$ a $C^1$ mapping 
from $J \times \mathcal{X}$ to $\mathcal{Y}$.
Assume that 
\begin{gather*}
D_{\lam w}^2 \mathcal{F} \mbox{ and } D_{ww}^2 \mathcal{F} 
\mbox{ exist and continuous in } J \times \mathcal{X},
\\
\mathcal{F}(\lam,0)=0 \quad \mbox{for all } \lam \in J,
\\
\dim K(D_w \mathcal{F}(\lam_0,0))=\codim R(D_w \mathcal{F}(\lam_0,0))=1
 \quad \mbox{for some } \lam_0 \in J,
\\
D_{\lam w}^2 \mathcal{F}(\lam_0,0) \varphi_0 \notin R(D_w \mathcal{F}(\lam_0,0)),
 \quad \mbox{where } K(D_w \mathcal{F}(\lam_0,0))=\operatorname{span}\{ \varphi_0\}.
\end{gather*}
Then there exist an open interval $\tilde J$ containing $0$,
a $C^1$ curve $\{ (\Lam (s),W(s))\}_{s \in \tilde J} \subset J \times \mathcal{X}$ 
and a neighborhood $\mathcal{N} \subset J \times \mathcal{X}$ of $(\lam_0,0)$ 
such that
\begin{equation}\label{lwpr}
\left\{
\begin{gathered}
(\Lam (0),W(0))=(\lam_0,0), \quad W_s(0)=\varphi_0,
\\
\mathcal{F}(\Lam (s),W(s))=0 \quad \mbox{for all } s \in \tilde J,
\\
\{ (\lam,w) \in \mathcal{N}; \mathcal{F}(\lam,w)=0\}
 =\left( \{ (\lam,0)\}_{\lam \in J} \cup \{ (\Lam (s),W(s))\}_{s \in \tilde J} \right) \cap {\cN},
\end{gathered}
\right.
\end{equation}
where $W_s$ stands for the derivative of $W$ with respect to $s$.
\end{thm}

We will use the following lemma to ensure 
the differentiability of an eigenpair of \eqref{llevp} with respect to $\lam$.
For the proof, 
we refer the reader to \cite[Proposition~I.7.2]{K12}.
\begin{lem}[\cite{K12}]\label{lem:devf}
In addition to the assumptions of Theorem~\ref{thm:crbt},
suppose that $\cX$ is continuously embedded in $\cY$
and that $R(D_w \mathcal{F}(\lam_0,0))$ is a complement of 
$K(D_w \mathcal{F}(\lam_0,0))=\operatorname{span}\{ \varphi_0\}$ in ${\cY}$:
\begin{equation*}
{\cY}=R(D_w \mathcal{F}(\lam_0,0)) \oplus \operatorname{span}\{ \varphi_0\}.
\end{equation*}
Then there exist an open interval $\tilde J$ containing $\lam_0$ and $C^1$ mappings
$\tilde J \ni \lam \mapsto \mu(\lam) \in \bbr$
and $\tilde J \ni \lam \mapsto \varphi (\lam) \in {\cX}$ such that 
\begin{equation*}
\mu(\lam_0)=0,
\qquad
\varphi (\lam_0)=\varphi_0,
\qquad
D_w \mathcal{F}(\lam,0) \varphi (\lam) =\mu(\lam) \varphi (\lam).
\end{equation*}
\end{lem}

\begin{remk}\label{rem:dmu}
We will apply Lemma~\ref{lem:devf} in the special case where
$\cY$ is also continuously embedded in some Hilbert space $\cZ$ and 
\begin{equation*}
R(D_w \mathcal{F}(\lam_0,0)) =\operatorname{span}\{ \varphi_0\}^{\perp} \cap {\cY}
 =\{ \varphi \in {\cY}; \langle \varphi,\varphi_0 \rangle=0 \}.
\end{equation*}
Here ${}^\perp$ and $\langle \cdot,\cdot \rangle$ 
stand for the orthogonal complement and the inner product in $\cZ$, respectively.
Then by differentiating the equality 
$\langle D_w \mathcal{F}(\lam,0) \varphi (\lam),\varphi_0 \rangle
=\mu(\lam) \langle \varphi (\lam),\varphi_0 \rangle$,
we obtain the well-known formula
\begin{equation}\label{fode}
\frac{d\mu}{d\lam}(\lam_0)
 =\frac{\langle D_{\lam w}^2 \mathcal{F}(\lam_0,0) \varphi_0,\varphi_0 \rangle}{
  \langle \varphi_0,\varphi_0 \rangle}.
\end{equation}
\end{remk}

To apply Theorem~\ref{thm:crbt} to our problem,
we set up function spaces ${\cX}$, ${\cY}$ and ${\cZ}$.
We choose
\begin{gather*}
{\cX}:=\left\{ u \in X_0 ; \,
\begin{aligned}
&u_x(-1)=u_x(1)=0
\\
&u(-0) +au_x(-0)=u(+0)-au_x(+0)
\\
&u_x(-0)=u_x(+0)
\end{aligned}
\right\},
\\
{\cY}:=\{ u \in C([-1,1] \setminus \{0\});
 u|_{[-1,0)} \mbox{ and } u|_{(0,1]} \mbox{ are uniformly continuous}\},
\\
{\cZ}:=L^2(-1,1).
\end{gather*}
Then ${\cX}$ is a closed linear subspace of $X_0$
and ${\cY}$ is a Banach space endowed with the uniform norm.
Let $\langle \cdot, \cdot \rangle$ denote the inner product on ${\cZ}$:
\begin{equation*}
\langle u,v \rangle:=\int_{-1}^1 u(x)v(x) dx,
\quad u,v \in {\cZ}.
\end{equation*}
We note that integrating by parts yields
\begin{align}
\langle \varphi_{xx}+q\varphi, \psi \rangle -\langle \varphi, \psi_{xx}+q\psi \rangle
 =&\big[\varphi_x \psi -\varphi \psi_x\big]^{x=-0}_{x=-1} 
  +\big[\varphi_x \psi -\varphi \psi_x\big]_{x=+0}^{x=1}
\nonumber
\\
 =&\big[\varphi_x \psi -\varphi \psi_x\big]_{x=-1}^{x=1}
  +\big\{ \varphi_x (\psi +a\psi_x) -(\varphi +a\varphi_x) \psi_x \big\} \big|_{x=-0}
\nonumber
\\
&-\big\{ \varphi_x (\psi -a\psi_x) -(\varphi -a\varphi_x) \psi_x \big\} \big|_{x=+0}
\label{Tsym}
\end{align}
for all $\varphi,\psi \in X_0$ and $q \in {\cY}$.
We define a $C^2$ mapping ${\cG}:(0,\infty) \times {\cX} \to {\cY}$ by
\begin{equation*}
{\cG} (\lam,u):=u_{xx} +\lam f(u).
\end{equation*}
By definition, we have 
${\cS}=\{ (\lam,u) \in (0,\infty) \times {\cX};{\cG} (\lam,u)=0, |u|<1\}$.

The following lemma can be shown in a standard way.
We give a proof in Appendix~\ref{appendixB} for readers' convenience.
\begin{lem}\label{lem:bt1}
Let $q \in {\cY}$ and define a linear operator $T:{\cX} \to {\cY}$ 
by $T:=d^2/dx^2 +q$.
Then
\begin{equation}\label{ReKp}
R(T)=K(T)^\perp \cap {\cY}
 =\{ \varphi \in {\cY}; \langle \varphi,\psi \rangle=0 \mbox{ for all } \psi \in K(T)\}.
\end{equation}
\end{lem}

Let us prove Proposition~\ref{prop:bt}.
\begin{proof}[Proof of Proposition~\ref{prop:bt}]
Set
\begin{equation*}
T_\lam:=D_u \mathcal{G}(\lam,u(\cdot,\lam))=\frac{d^2}{dx^2}+\lam f'(u(\cdot,\lam)).
\end{equation*}
From Lemma~\ref{lem:ndcD},
we see that the condition $D(\lam_0,\be_1(\lam_0),\be_2(\lam_0)) \neq 0$ 
gives $K(T_{\lam_0})=\{0\}$.
This with Lemma~\ref{lem:bt1} shows that 
the linear operator $T_{\lam_0}:{\cX} \to {\cY}$ 
is an isomorphism if $D(\lam_0,\be_1(\lam_0),\be_2(\lam_0)) \neq 0$.
Hence the assertion (i) follows from the implicit function theorem.

We prove (ii) and (iii).
In what follows, let \eqref{assmp:bt} hold.
We put
\begin{equation*}
\mathcal{F}(\lam,w):=\mathcal{G}(\lam,u(\cdot,\lam)+w),
 \quad (\lam,w) \in J \times {\cX},
\end{equation*}
and check that the conditions of Theorem~\ref{thm:crbt} and Lemma~\ref{lem:devf} hold.
We infer that ${\cF}:J \times {\cX} \to {\cY}$ is continuously differentiable
and has continuous derivatives $D^2_{\lam w} {\cF}$ and $D_{ww}^2 {\cF}$,
since $\cG$ is of class $C^2$ and 
the mapping $J \ni \lam \mapsto u(\cdot,\lam) \in {\cX}$ is of class $C^1$.
Moreover, by assumption, we have ${\cF}(\lam,0)=0$ for $\lam \in J$.
Since the condition $D(\lam_0,\be_1 (\lam_0),\be_2 (\lam_0))=0$ is assumed,
we see from Lemmas~\ref{lem:ndcD} and \ref{lem:ipmi} that 
for some $\varphi_0 \in {\cX} \setminus \{0\}$,
\begin{equation}\label{ktlz}
K(T_{\lam_0})=\operatorname{span}\{ \varphi_0\}.
\end{equation}
For $j=1,2$, let $\phi_j=\phi_j(x,\lam)$ be given by \eqref{phjdef} with $\be_j=\be_j(\lam)$.
We put $\al_1:=\varphi_0(-1)$, $\al_2:=\varphi_0(1)$ and
\begin{equation*}
\varphi (x,\lam):=\left\{
\begin{aligned}
&\al_1 \phi_1(x,\lam) && (x \in [-1,0)), \\
&\al_2 \phi_2(x,\lam) && (x \in (0,1]).
\end{aligned}
\right.
\end{equation*}
As shown in the proof of Lemma~\ref{lem:ndcD} (see \eqref{pjvpj}),
both $\al_1$ and $\al_2$ are nonzero and
the assumption $D(\lam_0,\be_1 (\lam_0),\be_2 (\lam_0))=0$ gives
\begin{equation}\label{pseph}
\varphi (\cdot,\lam_0)=\varphi_0.
\end{equation}
By \eqref{PQphj}, we have
\begin{equation*}
-\frac{\al_1 \al_2 \sqrt{\lam}}{G'(\be_1(\lam)) G'(\be_2(\lam))} D(\lam,\be_1(\lam),\be_2(\lam))
 =\det 
\begin{pmatrix}
(\varphi +a\varphi_x)|_{x=-0} & (\varphi -a\varphi_x)|_{x=+0} \\
\varphi_x|_{x=-0} & \varphi_x|_{x=+0}
\end{pmatrix}.
\end{equation*}
Differentiating this with respect to $\lam$,
we find that
\begin{align}
&-\frac{\al_1 \al_2 \sqrt{\lam_0}}{G'(\be_1(\lam_0)) G'(\be_2(\lam_0))} 
 \left. \frac{d}{d\lam} D(\lam,\be_1(\lam),\be_2(\lam)) \right|_{\lam=\lam_0}
\nonumber
\\
&=\left. \det 
\begin{pmatrix}
(\varphi_\lam +a\varphi_{x\lam})_{x=-0} & (\varphi_\lam -a\varphi_{x\lam})|_{x=+0} \\
(\varphi_0)_x|_{x=-0} & (\varphi_0)_x|_{x=+0}
\end{pmatrix}
\right|_{\lam=\lam_0}
\nonumber
\\
&\quad +\left. \det 
\begin{pmatrix}
\left. \left\{ \varphi_0 +a(\varphi_0)_x \right\} \right|_{x=-0}
 & \left. \left\{ \varphi_0 -a(\varphi_0)_x \right\} \right|_{x=+0} \\
\varphi_{x\lam} |_{x=-0} & \varphi_{x\lam} |_{x=+0}
\end{pmatrix}
\right|_{\lam=\lam_0},
\label{dDps}
\end{align}
where we have used the assumption $D(\lam_0,\be_1 (\lam_0),\be_2 (\lam_0))=0$
and \eqref{pseph}.

Since $\phi_j$ satisfies \eqref{duaeq},
we see that $\varphi_{xx}+\lam f'(u(\cdot,\lam))\varphi =0$ in $(-1,1) \setminus \{0\}$
and $\varphi_x=0$ at $x=\pm 1$. 
Hence, using \eqref{Tsym} 
for $\varphi=\varphi_0$, $\psi=\varphi(\cdot,\lam)$ and $q=\lam f'(u(\cdot,\lam))$,
we have
\begin{align*}
\langle T_\lam \varphi_0, \varphi \rangle
 =&\big[ (\varphi_0)_x (\varphi +a\varphi_x)
  -\{ \varphi_0 +a(\varphi_0)_x\} \varphi_x \big] \big|_{x=-0}
\\
&-\big[ (\varphi_0)_x (\varphi -a\varphi_x)
 -\{ \varphi_0 -a(\varphi_0)_x\} \varphi_x \big] \big|_{x=+0}.
\end{align*}
We differentiate this equality with respect to $\lam$.
By \eqref{ktlz}, the derivative of the left-hand side at $\lam=\lam_0$ is computed as 
\begin{equation*}
\left. \frac{d}{d\lam} \langle T_\lam \varphi_0, \varphi \rangle \right|_{\lam=\lam_0}
 =\left. \frac{d}{d\lam} \langle T_\lam \varphi_0, \varphi_0 \rangle \right|_{\lam=\lam_0}
  +\left. \langle T_{\lam_0} \varphi_0, \varphi_\lam \rangle \right|_{\lam=\lam_0}
   =\langle D_{\lam w}^2 \mathcal{F}(\lam_0,0) \varphi_0,\varphi_0 \rangle.
\end{equation*}
Therefore
\begin{align*}
\langle D_{\lam w}^2 \mathcal{F}(\lam_0,0) \varphi_0,\varphi_0 \rangle
 =&\big[ (\varphi_0)_x (\varphi_\lam +a\varphi_{x\lam})
  -\{ \varphi_0 +a(\varphi_0)_x\} \varphi_{x\lam} \big] \big|_{x=-0,\lam=\lam_0}
\\
&-\big[ (\varphi_0)_x (\varphi_\lam -a\varphi_{x\lam})
 -\{ \varphi_0 -a(\varphi_0)_x\} \varphi_{x\lam} \big] \big|_{x=+0,\lam=\lam_0}.
\end{align*}
Notice that the right-hand side of this equality coincides with that of \eqref{dDps}, 
since $\varphi_0$ satisfies 
\begin{equation*}
\{ \varphi_0 +a(\varphi_0)_x\} |_{x=-0}=\{ \varphi_0 -a(\varphi_0)_x\} |_{x=+0},
\quad
(\varphi_0)_x|_{x=-0}=(\varphi_0)_x|_{x=+0}.
\end{equation*}
Thus
\begin{equation}\label{rncd}
\langle D_{\lam w}^2 \mathcal{F}(\lam_0,0) \varphi_0,\varphi_0 \rangle
 =-\frac{\al_1 \al_2 \sqrt{\lam_0}}{G'(\be_1(\lam_0)) G'(\be_2(\lam_0))}
  \left. \frac{d}{d\lam} D(\lam,\be_1(\lam),\be_2(\lam)) \right|_{\lam=\lam_0}.
\end{equation}

From \eqref{rncd} and the second condition of \eqref{assmp:bt}, 
we see that $\langle D_{\lam w}^2 \mathcal{F}(\lam_0,0) \varphi_0,\varphi_0 \rangle \neq 0$.
Combining this with Lemma~\ref{lem:bt1} and \eqref{ktlz} gives
\begin{equation*}
D_{\lam w}^2 \mathcal{F}(\lam_0,0) \varphi_0
 \notin \operatorname{span}\{ \varphi_0\}^{\perp} \cap {\cY}
  =R(T_{\lam_0})=R(D_w \mathcal{F}(\lam_0,0)).
\end{equation*}
Hence all the conditions of Theorem~\ref{thm:crbt} and Lemma~\ref{lem:devf} are satisfied.
Applying Theorem~\ref{thm:crbt},
we obtain a $C^1$ curve 
$\{ (\Lam(s),W(\cdot,s))\}_{s \in \tilde J} \subset J \times {\cX}$ satisfying \eqref{lwpr}.
Then one can directly check that 
$\tilde {\cC}:=\{ (\Lam(s),u(\cdot,\Lam(s))+W(\cdot,s))\}_{s \in \tilde J}$
is the one having the desired properties stated in (ii).

It remains to derive \eqref{evdf}.
By Lemma~\ref{lem:devf}, 
we have an eigenvalue $\mu(\lam)$ of $D_w {\cF}(\lam,0)=T_\lam$
which is of class $C^1$ and satisfies $\mu(\lam_0)=0$.
It is easily seen that $\mu(\lam)$ is an eigenvalue of \eqref{llevp} for $u=u(\cdot,\lam)$.
Hence $\mu (\lam)$ coincides with $\mu_n (u(\cdot,\lam))$, 
since Lemma~\ref{lem:ipmi} shows that 
each eigenvalue $\mu_m (u(\cdot,\lam))$, $m \in \bbn \cup \{0\}$, 
is isolated and continuous with respect to $\lam$.
In particular, $\mu_n (u(\cdot,\lam))$ is continuously differentiable 
in a neighborhood of $\lam_0$.
As noted in Remark~\ref{rem:dmu},
we can compute the derivative of $\mu (\lam)=\mu_n (u(\cdot,\lam))$ 
by the formula \eqref{fode}.
Therefore we see from \eqref{rncd} and Lemma~\ref{lem:veef} that
\begin{align*}
\sgn \left( \left. \frac{d}{d\lam} \mu_n (u(\cdot,\lam)) \right|_{\lam=\lam_0} \right)
 &=\sgn \left( -\varphi_0 (-1)\varphi_0 (1)
  \left. \frac{d}{d\lam} D(\lam,\be_1(\lam),\be_2(\lam)) \right|_{\lam=\lam_0} \right)
\\
&=(-1)^{n-1} \sgn \left( 
 \left. \frac{d}{d\lam} D(\lam,\be_1(\lam),\be_2(\lam)) \right|_{\lam=\lam_0} \right).
\end{align*}
We thus obtain \eqref{evdf}, and the proof is complete.
\end{proof}


\section{Primary branches}\label{sec:pb}

In this section, we examine primary branches of solutions of \eqref{leq}
bifurcating from the trivial branch ${\cL}=\{ (\lam,0)\}_{\lam \in (0,\infty)}$.
To set up notation, we introduce the function
\begin{equation*}
g(\be,\phi) :=\left\{
\begin{aligned}
&\frac{G(\be \cos \phi)}{\be}-a\sin \phi \int_0^{\phi}G'(\be \cos \tau)d\tau
&& \mbox{if } \be \neq 0,
\\
&G'(0)\cos \phi -aG'(0) \phi \sin \phi
&& \mbox{if } \be=0.
\end{aligned}
\right.
\end{equation*}
In Lemma~\ref{lem:sgez} below, we will show that
for $k \in \bbn$, there exists $\phi_k=\phi_k(\be) \in C^1(I)$ satisfying
$\phi_k(\be) \in ((k-1)\pi,(k-1/2)\pi)$ and $g(\be,\phi_k(\be))=0$.
Then we put
\begin{gather*}
\lam_k^o(\be):=\left( \int_0^{\phi_k (\be)} G'(\be \cos \tau) d\tau \right)^2,
\qquad
\lam_k^e(\be):=\left( \int_0^{k\pi} G'(\be \cos \tau) d\tau \right)^2,
\\
u^o_k(x.\be):=\Psi_{\lam_k^o(\be)} (\be,-\be)(x),
\qquad
u^e_k(x,\be):=\Psi_{\lam_k^e(\be)} (\be,\be)(x),
\\
{\cC}_k^o:=\{ (\lam_k^o (\be),u^o_k(\cdot,\be))\}_{\be \in I},
\qquad
{\cC}_k^e:=\{ (\lam_k^e (\be),u^e_k(\cdot,\be))\}_{\be \in I}.
\end{gather*}
By the definitions of $g$ and $z_k$, we have 
\begin{equation}\label{phpr1}
\phi_k (0)=z_k.
\end{equation}
This together with the fat that $G'(0)=1/\sqrt{f'(0)}$ (see \eqref{Gdao} below) gives
\begin{gather*}
\lam_k^o(0)=(\phi_k (0) G'(0))^2=\lam_{2k-1},
\qquad
\lam_k^e(0):=(k\pi G'(0))^2=\lam_{2k}.
\end{gather*}
Therefore ${\cC}_k^o$ (resp. ${\cC}_k^e$) is a $C^1$ curve in $(0,\infty) \times X$
which intersects $\cL$ at $(\lam_k^o(0),u^o_k(\cdot,0))=(\lam_{2k-1},0)$ 
(resp. $(\lam_k^e(0),u^e_k(\cdot,0))=(\lam_{2k},0)$).
We define $C_{k,+}^o$, $C_{k,-}^o$, $C_{k,+}^e$ and $C_{k,-}^e$ by
\begin{gather*}
{\cC}_{k,+}^o:=\{ (\lam_k^o (\be),u^o_k(\cdot,\be))\}_{\be \in (0,\be_0)},
\qquad
{\cC}_{k,-}^o:=\{ (\lam_k^o (\be),u^o_k(\cdot,\be))\}_{\be \in (-\be_0,0)},
\\
{\cC}_{k,+}^e:=\{ (\lam_k^e (\be),u^e_k(\cdot,\be))\}_{\be \in (0,\be_0)},
\qquad
{\cC}_{k,-}^e:=\{ (\lam_k^e (\be),u^e_k(\cdot,\be))\}_{\be \in (-\be_0,0)}.
\end{gather*}

We prove the following proposition.
\begin{prop}\label{prop:pb}
The following hold.
\begin{enumerate}[(i)]
\item
There hold
\begin{gather}
{\cC}_{k,-}^o=-{\cC}_{k,+}^o,
\qquad
{\cC}_{k,-}^e=-{\cC}_{k,+}^e,
\label{cksy}
\\
\ \bigcup_{k=1}^\infty {\cC}_{k,+}^o \cup {\cC}_{k,-}^o =\mathcal{S}^o,
\qquad 
\bigcup_{k=1}^\infty {\cC}_{k,+}^e \cup {\cC}_{k,-}^e =\mathcal{S}^e.
\label{cues}
\end{gather}
Furthermore, for every $\lam \in (0,\infty)$, there exists a neighborhood 
$\cN$ of $(\lam,0)$ in $(0,\infty) \times X_0$ such that
\begin{equation}\label{lucl}
{\cS} \cap {\cN}=\left\{
\begin{aligned}
&({\cC}_k^o \cup {\cL}) \cap {\cN} && \mbox{if } \lam=\lam_{2k-1},
\\
&({\cC}_k^e \cup {\cL}) \cap {\cN} && \mbox{if } \lam=\lam_{2k},
\\
&{\cL} \cap {\cN} && \mbox{if } \lam \neq \lam_{2k-1}, \lam_{2k}.
\end{aligned}
\right.
\end{equation}
\item
Assume \eqref{aasfw}. Then
\begin{equation}\label{sdlk}
\sgn (\be) \frac{d\lam_k^o}{d\be} (\be)>0,
\quad
\sgn (\be) \frac{d\lam_k^e}{d\be} (\be)>0
\quad \mbox{for } \be \in I \setminus \{0\}.
\end{equation}
In particular, $C_{k,+}^o$ and $C_{k,-}^o$ (resp. $C_{k,+}^e$ and $C_{k,-}^e$)
are parametrized by $\lam \in (\lam_{2k-1},\infty)$ (resp. $\lam \in (\lam_{2k},\infty)$).
\end{enumerate}
\end{prop}


We begin with solving the equation $g(\be,\phi)=0$.
\begin{lem}\label{lem:sgez}
For $k \in \bbn$ and $\be \in I$,
there exists $\phi_k (\be) \in ((k-1)\pi,(k-1/2)\pi)$ with the following properties:
\begin{gather}
\{ \phi \in [0,\infty); g(\be,\phi)=0\} =\{ \phi_k(\be))\}_{k=1}^\infty,
\label{phpr3}
\\
\phi_k (-\be)=\phi_k (\be),
\label{phpr4}
\\
\lim_{\be \to \pm \be_0} \phi_k (\be)=(k-1)\pi.
\label{phpr2}
\end{gather}
Furthermore, $\phi_k \in C^1(I)$ and
\begin{equation}\label{dphk}
\frac{d\phi_k}{d\be}(\be) =\frac{J(\be,\phi_k(\be))}{I(\be,\phi_k(\be))},
\end{equation}
where
\begin{gather*}
I(\be,\phi):=\be \left\{ \left( \frac{G'(\be \cos \phi) \be \sin \phi}{G(\be \cos \phi)}
 +\frac{\cos \phi}{\sin \phi} \right) \int_0^{\phi}G'(\be \cos \tau)d\tau
  +G'(\be \cos \phi) \right\},
\\
J(\be,\phi):=\left( \frac{G'(\be \cos \phi) \be \cos \phi}{G(\be \cos \phi)} -1 \right)
 \int_0^\phi G'(\be \cos \tau)d\tau
  -\be \int_0^\phi G''(\be \cos \tau) \cos \tau d\tau.
\end{gather*}
\end{lem}

\begin{proof}
For the moment suppose that $k$ is odd.
We note that $g \in C^1(I \times \bbr)$,
since $g$ is written as
\begin{equation*}
g(\be,\phi)=\cos \phi \int_0^1 G'(t\be \cos \tau) dt
 -a\sin \phi \int_0^{\phi}G'(\be \cos \tau)d\tau.
\end{equation*}
From the fact that $G'$ is positive, we have
\begin{equation*}
g(\be,(k-1)\pi)>0,
\qquad
g(\be,\phi)<0 \quad \mbox{for } \phi \in \left[ \left( k-\frac{1}{2}\right) \pi,k\pi \right].
\end{equation*}
Furthermore, 
\begin{align*}
g_\phi (\be,\phi) &=-(1+a)G'(\be \cos \phi) \sin \phi 
 -a\cos \phi \int_0^{\phi}G'(\be \cos \tau)d\tau 
\\
&<0 \quad \mbox{for } \phi \in \left( (k-1) \pi,\left( k-\frac{1}{2}\right)\pi \right).
\end{align*}
Therefore there exists $\phi_k(\be) \in ((k-1)\pi,(k-1/2) \pi)$
such that 
\begin{equation*}
\{ \phi \in [(k-1)\pi,k\pi]; g(\be,\phi)=0\}=\{ \phi_k(\be)\},
\end{equation*}
and the implicit function theorem shows that $\phi_k \in C^1(I)$.
The case where $k$ is even can be dealt with in the same way.
We have thus shown the existence of $\phi_k(\be)$ and \eqref{phpr3}.

The fact that $G$ is odd yields $g(-\be,\phi)=g(\be,\phi)$,
and hence $g(-\be,\phi_k(\be))=g(\be,\phi_k(\be))=0$.
From \eqref{phpr3}, we obtain \eqref{phpr4}.

We prove \eqref{phpr2}.
To obtain a contradiction,
suppose that $\phi_k (\be)$ does not converge to $(k-1)\pi$ 
as $\be \to \be_0$ or $\be \to -\be_0$.
Then we can take $\{ \be_j\}_{j=1}^\infty \subset I$ 
and $\del \in (0,\pi/2)$ such that
\begin{equation*}
|\be_j|<|\be_{j+1}|,
\qquad
|\be_j| \to \be_0 \quad (j \to \infty),
\qquad 
\phi_k(\be_j) \ge (k-1)\pi +\del.
\end{equation*}
By \eqref{Gpas} and the monotone convergence theorem, 
\begin{align*}
\int_0^{\phi_k(\be_j)}G'(\be_j \cos \tau)d\tau
 &\ge c\int_0^{(k-1)\pi +\del} \frac{d\tau}{\sqrt{\be_0-|\be_j| \cos \tau}}
\\
&\to c\int_0^{(k-1)\pi +\del} \frac{d\tau}{\sqrt{\be_0-\be_0 \cos \tau}}
 =\infty \quad (j \to \infty).
\end{align*}
From this and the fact that $|\sin \phi_k(\be_j)| \ge |\sin \del|$,
we have $|g(\be_j,\phi_k(\be_j))| \to \infty$ as $j \to \infty$. 
This contradicts the equality $g(\be_j,\phi_k(\be_j))=0$,
and therefore \eqref{phpr2} holds.

It remains to prove \eqref{dphk}.
Differentiating $g(\be,\phi_k(\be))=0$ yields
\begin{align}
\frac{d\phi_k}{d\be}(\be) &=-\frac{g_\be (\be,\phi_k(\be))}{g_\phi (\be,\phi_k(\be))}
\nonumber
\\
&=\left. \frac{G'(\be \cos \phi) \cos \phi -G(\be \cos \phi)/\be
 -a\be \sin \phi \int_0^{\phi}G''(\be \cos \tau) \cos \tau d\tau}{
  (1+a)G'(\be \cos \phi) \be \sin \phi +a\be \cos \phi \int_0^{\phi}G'(\be \cos \tau)d\tau}
   \right|_{\phi=\phi_k(\be)}.
\label{dphk0}
\end{align}
We note that
\begin{equation}\label{ge0e}
\frac{1}{a}=\left. \frac{\be \sin \phi}{G(\be \cos \phi)}
 \int_0^{\phi}G'(\be \cos \tau)d\tau \right|_{\phi=\phi_k(\be)},
\end{equation}
which follows from $g(\be,\phi_k(\be))=0$.
Plugging this into \eqref{dphk0}, 
we obtain \eqref{dphk}.
Therefore the lemma follows.
\end{proof}

We prove the further property of $\phi_k$ to be used later.
\begin{lem}\label{lem:phli}
For each $k \in \bbn$,
\begin{equation*}
\liminf_{\be \to \pm \be_0} \frac{|\sin \phi_k (\be)|}{\sqrt{\be_0-|\be|}}>0.
\end{equation*}
\end{lem}

\begin{proof}
We use \eqref{Gpas} to obtain
\begin{equation*}
G'(\be \cos \tau) \le \frac{C}{\sqrt{\be_0-|\be \cos \tau|}}
 \le \frac{C}{\sqrt{\be_0-|\be|}}.
\end{equation*}
Hence, by \eqref{ge0e}, we have
\begin{equation*}
\frac{1}{a} \le \frac{C|\be| \phi_k(\be)}{|G(\be \cos \phi_k(\be))|}
 \cdot \frac{|\sin \phi_k (\be)|}{\sqrt{\be_0-|\be|}}.
\end{equation*}
The desired inequality is verified by combining this with \eqref{phpr2}.
\end{proof}

%

It is known (see \cite{O61,S90,SW81}) that the condition \eqref{aasfw} implies the inequality
\begin{equation*}
\sgn(\be) \int_0^{l\pi} G''(\be \cos \tau) \cos \tau d\tau>0
 \quad (\be \in I \setminus \{0\}, l \in \bbn),
\end{equation*}
which is equivalent to
\begin{equation*}
\sgn(\al) \frac{d}{d\al} \int_0^\al \frac{ds}{\sqrt{F(\al)-F(s)}}>0
 \quad (\al \in (-1,1) \setminus \{0\}).
\end{equation*}
In order to show \eqref{sdlk},
we generalize this inequality. 
\begin{lem}\label{lem:gdpe}
Let \eqref{aasfw} hold.
Then for all $\be \in I \setminus \{0\}$ and $\phi \in (0,\infty)$,
\begin{align*}
&\be \int_0^\phi G''(\be \cos \tau) \cos \tau d\tau 
\\
&>\left\{
\begin{aligned}
&-\left( G'(\be \cos \phi) \cos \phi
 -\frac{G'(\be)}{G(\be)} G(\be \cos \phi) \right) \frac{1}{\sin \phi}
  && \mbox{if } \sin \phi \neq 0,
\\
&0 && \mbox{if } \sin \phi=0.
\end{aligned}
\right.
\end{align*}
\end{lem}

\begin{proof}
We first consider the case $\sin \phi \neq 0$.
Put
\begin{equation*}
\tilde H(v):=\left\{
\begin{aligned}
&\frac{G(v)}{v} && (v \neq 0),
\\
&G'(0) && (v=0).
\end{aligned}
\right.
\end{equation*}
From Lemma~\ref{lem:Gpro}, 
one can check that $\tilde H \in C^1(I)$ and $\tilde H'(v)=G'(v)H(v)/v^3$.
Then
\begin{align}
\be \int_0^\phi G''(\be \cos \tau) \cos \tau d\tau
 &=\be \int_0^\phi \left( G''(\be \cos \tau)
  -\frac{G'(\be) \be}{G(\be)} \tilde H'(\be \cos \tau) \right) \cos \tau d\tau
\nonumber
\\
&\quad +\frac{G'(\be)\be^2}{G(\be)} \int_0^\phi \tilde H'(\be \cos \tau) \cos \tau d\tau
\nonumber
\\
&=-\int_0^\phi \frac{\D}{\D \tau} \left( G'(\be \cos \tau)
 -\frac{G'(\be) \be}{G(\be)} \tilde H(\be \cos \tau) \right)
  \cdot \frac{\cos \tau}{\sin \tau} d\tau
\nonumber
\\
&\quad +\frac{G'(\be)}{G(\be)\be}
 \int_0^\phi \frac{H(\be \cos \tau)G'(\be \cos \tau)}{\cos^2 \tau} d\tau.
\label{gdpe1}
\end{align}
We apply integration by parts to obtain
\begin{align}
&\int_0^\phi \frac{\D}{\D \tau} \left( G'(\be \cos \tau)
 -\frac{G'(\be) \be}{G(\be)} \tilde H(\be \cos \tau) \right)
  \cdot \frac{\cos \tau}{\sin \tau} d\tau
\nonumber
\\
&=\left( G'(\be \cos \phi)
 -\frac{G'(\be) \be}{G(\be)} \tilde H(\be \cos \phi) \right) \frac{\cos \phi}{\sin \phi}
\nonumber
\\
&\quad +\int_0^\phi \left( G'(\be \cos \tau)
 -\frac{G'(\be) \be}{G(\be)} \tilde H(\be \cos \tau) \right) \frac{1}{\sin^2 \tau} d\tau.
\label{gdpe2}
\end{align}
This computation is valid, since
\begin{equation}\label{ghas}
G'(\be \cos \phi)
 -\frac{G'(\be) \be}{G(\be)} \tilde H(\be \cos \phi) =O(1-|\cos \tau|)=O(\sin^2 \tau)
\quad
\mbox{as } \tau \to l\pi, \ l \in \bbz.
\end{equation}
Note that the integrand of the second term on the right of \eqref{gdpe2}
is written as
\begin{align*}
G'(\be \cos \tau) -\frac{G'(\be) \be}{G(\be)} \tilde H(\be \cos \tau)
&=\frac{G'(\be)G'(\be \cos \tau)}{G(\be) \be}
 \left( \frac{G(\be) \be}{G'(\be)} -
  \frac{G(\be \cos \tau) \be}{G'(\be \cos \tau) \cos \tau} \right)
\\
&=\frac{G'(\be)G'(\be \cos \tau)}{G(\be) \be}
 \left( \frac{H(\be \cos \tau)}{\cos^2 \tau} -H(\be) \right).
\end{align*}
Therefore
\begin{align*}
&\int_0^\phi \frac{\D}{\D \tau} \left( G'(\be \cos \tau)
 -\frac{G'(\be) \be}{G(\be)} \tilde H(\be \cos \tau) \right)
  \cdot \frac{\cos \tau}{\sin \tau} d\tau
\\
&=\left( G'(\be \cos \phi) \cos \phi
 -\frac{G'(\be)}{G(\be)} G(\be \cos \phi) \right) \frac{1}{\sin \phi}
\\
&\quad +\frac{G'(\be)}{G(\be) \be} 
 \int_0^\phi \left( \frac{H(\be \cos \tau)}{\cos^2 \tau} -H(\be) \right)
  \frac{G'(\be \cos \tau)}{\sin^2 \tau} d\tau.
\end{align*}
Substituting this into \eqref{gdpe1},
we find that
\begin{align*}
\be \int_0^\phi G''(\be \cos \tau) \cos \tau d\tau
 &=-\left( G'(\be \cos \phi) \cos \phi
  -\frac{G'(\be)}{G(\be)} G(\be \cos \phi) \right) \frac{1}{\sin \phi}
\\
&\quad +\frac{G'(\be)}{G(\be)\be} \int_0^\phi
 \frac{G'(\be \cos \tau)(H(\be) -H(\be \cos \tau))}{\sin^2 \tau} d\tau.
\end{align*}
Since the assumption \eqref{aasfw} implies \eqref{aasGw},
we deduce that the second term on the right of this equality is positive.
Thus we obtain the desired inequality.

In the case $\sin \phi =0$,
we see from \eqref{ghas} that the first term on the right of \eqref{gdpe2} vanishes,
and hence the same argument works.
\end{proof}

To obtain odd and even solutions of \eqref{leq},
we find solutions of \eqref{smteq} satisfying either $\be_1=-\be_2$ or $\be_1=\be_2$.
\begin{lem}\label{lem:PQzs}
There hold
\begin{gather*}
\{ (\lam,\be) \in (0,\infty) \times (I \setminus \{0\}); (\lam,\be,-\be) \in {\cT} \} 
 =\bigcup_{k=1}^\infty \{ (\lam_k^o(\be),\be) \}_{\be \in I \setminus \{0\}},
\\
\{ (\lam,\be) \in (0,\infty) \times (I \setminus \{0\}); (\lam,\be,\be) \in {\cT} \} 
 =\bigcup_{k=1}^\infty \{ (\lam_k^e(\be),\be) \}_{\be \in I \setminus \{0\}}.
\end{gather*}
\end{lem}

\begin{proof}
We see from \eqref{PQsy} that 
$(\lam,\be,-\be) \in {\cT}$ if and only if $P(\lam,\be)=0$,
and that $(\lam,\be,\be) \in {\cT}$ if and only if $Q(\lam,\be)=0$.
Hence
\begin{gather*}
\{ (\lam,\be); (\lam,\be,-\be) \in {\cT}, \be \neq 0\} 
 =\{ (\lam,\be); P(\lam,\be)=0, \be \neq 0\},
\\
\{ (\lam,\be); (\lam,\be,\be) \in {\cT}, \be \neq 0\} 
 =\{ (\lam,\be); Q(\lam,\be)=0, \be \neq 0\}.
\end{gather*}
We fix $\be \in I \setminus \{0\}$.
By definition, we have $P(\lam,\be)=\be g(\be,\theta(\lam,\be))$.
This together with \eqref{phpr3} and \eqref{tblr} yields
\begin{equation*}
\{ \lam ; P(\lam,\be)=0\}
 =\{ \lam ; g(\be,\theta(\lam,\be))=0\}
  =\bigcup_{k=1}^\infty \{ \lam ; \theta (\lam,\be)=\phi_k (\be)\}
   =\bigcup_{k=1}^\infty \{ \lam_k^o(\be)\}.
\end{equation*}
Moreover, by \eqref{tblr},
\begin{equation*}
\{ \lam ; Q(\lam,\be)=0\}
 =\{ \lam ; \sin \theta (\lam,\be)=0\}
  =\bigcup_{k=1}^\infty \{ \lam ; \theta (\lam,\be)=k\pi \}
   =\bigcup_{k=1}^\infty \{ \lam_k^e(\be)\}.
\end{equation*}
Therefore we obtain the desired conclusion.
\end{proof}

Put
\begin{equation*}
\tilde {\cS}_\lam^o:=\{ u \in {\cS}_\lam \setminus \{0\}; u(-1)=-u(1) \},
\qquad
\tilde {\cS}_\lam^e:=\{ u \in {\cS}_\lam \setminus \{0\}; u(-1)=u(1) \}.
\end{equation*}

\begin{lem}\label{lem:sose}
There hold ${\cS}_\lam^o=\tilde {\cS}_\lam^o$ and ${\cS}_\lam^e=\tilde {\cS}_\lam^e$.
\end{lem}

\begin{proof}
It is clear that ${\cS}_\lam^o \subset \tilde {\cS}_\lam^o$.
To show $\tilde {\cS}_\lam^o \subset {\cS}_\lam^o$,
we let $u \in \tilde {\cS}_\lam^o$.
Then $u_1:=u|_{[-1,0)}$ and $u_2:=u|_{(0,1]}$ satisfy \eqref{u12ivp} 
for $\be_1=G^{-1}(u(-1))$ and $\be_2=G^{-1}(u(1))$.
Since the assumption $u \in \tilde {\cS}_\lam^o$ yields $\be_1=-\be_2$, 
we see that $u_1(x)$ and $-u_2(-x)$ satisfy the same initial value problem.
Hence $u_1(x)=-u_2(-x)$, which gives $u \in {\cS}_\lam^o$.
We have thus proved ${\cS}_\lam^o=\tilde {\cS}_\lam^o$.
The equality ${\cS}_\lam^e=\tilde {\cS}_\lam^e$ can be shown in the same way.
\end{proof}

We are now in a position to prove Proposition~\ref{prop:pb}.
\begin{proof}[Proof of Proposition~\ref{prop:pb}]
By \eqref{Pssy}, \eqref{phpr4} and the fact that $G$ is odd, 
we have 
\begin{equation}\label{lksy}
(\lam_k^o (-\be),u^o_k(\cdot,-\be))=(\lam_k^o (\be),-u^o_k(\cdot,\be)),
\quad
(\lam_k^e (-\be),u^e_k(\cdot,-\be))=(\lam_k^e (\be),-u^e_k(\cdot,\be)).
\end{equation}
Hence \eqref{cksy} follows.

Lemmas~\ref{ppoo} and \ref{lem:sose} yield
\begin{gather*}
{\cS}_\lam^o=\tilde {\cS}_\lam^o=
 \{ \Psi_\lam (\be,-\be); (\be,-\be) \in \mathcal{T}_\lam, \be \neq 0 \},
\\
{\cS}_\lam^e=\tilde {\cS}_\lam^e=
 \{ \Psi_\lam (\be,\be); (\be,\be) \in \mathcal{T}_\lam, \be \neq 0 \}.
\end{gather*}
Combining these with Lemma~\ref{lem:PQzs} shows that
\begin{align*}
{\cS}^o=\{ (\lam,\Psi_\lam (\be,-\be)); (\lam,\be,-\be) \in \mathcal{T}, \be \neq 0\}
 =\bigcup_{k=1}^\infty \{ (\lam_k^o (\be),u^o_k(\cdot,\be))\}_{\be \in I \setminus \{0\}},
\\
{\cS}^e=\{ (\lam,\Psi_\lam (\be,\be)); (\lam,\be,\be) \in \mathcal{T}, \be \neq 0\}
 =\bigcup_{k=1}^\infty \{ (\lam_k^e (\be),u^e_k(\cdot,\be))\}_{\be \in I \setminus \{0\}}.
\end{align*}
Therefore we obtain \eqref{cues}.

By \eqref{tblr} and the fact that $G'(0)=1/\sqrt{f'(0)}$ (see \eqref{Gdao} below), 
we deduce that $\theta(\lam,0)=\sqrt{\lam}/G'(0)=\sqrt{f'(0)\lam}$.
Hence it follows from \eqref{pqdb} that
\begin{align*}
D(\lam,0,0)&=2P_\be(\lam,0)Q_\be (\lam,0)
\\
&=-\frac{2}{\sqrt{f'(0)}} \left( \cos \sqrt{f'(0)\lam}
 -a\sqrt{f'(0)\lam} \sin \sqrt{f'(0)\lam}  \right) \sin \sqrt{f'(0)\lam}.
\end{align*}
This shows that $D(\lam,0,0)=0$ if and only if
$\sqrt{f'(0)\lam}=z_k$ or $\sqrt{f'(0)\lam}=k\pi$ for some $k \in \bbn$.
Moreover, 
\begin{gather*}
\left. \frac{d}{d\lam}  D(\lam,0,0) \right|_{\lam=z_k^2/f'(0)}
 =\frac{\sqrt{f'(0)}}{z_k} \{ (1+a)\sin z_k +a z_k \cos z_k\} \sin z_k>0,
\\
\left. \frac{d}{d\lam}  D(\lam,0,0) \right|_{\lam=(k\pi)^2/f'(0)}
 =-\frac{\sqrt{f'(0)}}{k\pi} <0.
\end{gather*}
Thus, using (i) and (ii) of Proposition~\ref{prop:bt} for ${\cC}={\cL}$,
we conclude that \eqref{lucl} holds.

What is left is to show (ii).
Assume \eqref{aasfw} and let $\be \in I \setminus \{0\}$.
The estimate for $\lam_k^e(\be)$ is directly derived from Lemma~\ref{lem:gdpe}.
Indeed,
\begin{equation*}
\be \frac{d}{d\be} \sqrt{\lam_k^e(\be)}
 =\be \int_0^{k\pi} G''(\be \cos \tau)\cos \tau d\tau>0.
\end{equation*}
Let us consider the estimate for $\lam_k^o(\be)$.
We use \eqref{dphk} to obtain
\begin{align*}
\be \frac{d}{d\be} \sqrt{\lam_k^o(\be)}
 &=\be \left( \int_0^{\phi_k(\be)} G''(\be \cos \tau) \cos \tau d\tau
  +\frac{d\phi_k}{d\be}(\be) G'(\be \cos \phi_k (\be)) \right) \\
&=\frac{\be}{I(\be,\phi)} \left. \left( I(\be,\phi) \int_0^{\phi} G''(\be \cos \tau) \cos \tau d\tau
 +J(\be,\phi) G'(\be\cos \phi) \right) \right|_{\phi=\phi_k(\be)}.
\end{align*}
Let $\phi \in ((k-1)\pi,(k-1/2)\pi)$. 
Then
\begin{equation}\label{sgni}
\frac{I(\be,\phi)}{\be} =\left( \frac{G'(\be \cos \phi) \be \sin \phi}{G(\be \cos \phi)}
 +\frac{\cos \phi}{\sin \phi} \right) \int_0^{\phi}G'(\be \cos \tau)d\tau +G'(\be \cos \phi) >0.
\end{equation}
By a direct computation, we have
\begin{align}
&\left. \left( I(\be,\phi) \int_0^{\phi} G''(\be \cos \tau) \cos \tau d\tau
 +J(\be,\phi) G'(\be\cos \phi) \right) \right/ \int_0^{\phi}G'(\be \cos \tau)d\tau
\nonumber
\\
&=\left( \frac{G'(\be \cos \phi) \be \sin \phi}{G(\be \cos \phi)}
 +\frac{\cos \phi}{\sin \phi} \right) 
  \be \int_0^{\phi}G''(\be \cos \tau) \cos \tau d\tau
\nonumber
\\
&\quad +\left( \frac{G'(\be \cos \phi) \be \cos \phi}{G(\be \cos \phi)} -1 \right)
 G'(\be \cos \phi).
\label{igjg0}
\end{align}
Applying Lemma~\ref{lem:gdpe} shows that the right-hand side of this equality is bounded below by
\begin{align*}
&\left( \frac{G'(\be \cos \phi) \be \sin \phi}{G(\be \cos \phi)}
 +\frac{\cos \phi}{\sin \phi} \right) \cdot \left\{ -\left( G'(\be \cos \phi) \cos \phi
  -\frac{G'(\be)}{G(\be)} G(\be \cos \phi) \right) \frac{1}{\sin \phi} \right\}
\\
&\quad +\left( \frac{G'(\be \cos \phi) \be \cos \phi}{G(\be \cos \phi)} -1 \right)
 G'(\be \cos \phi)
\\
&=\left( \be^2 \sin^2 \phi +\frac{G(\be \cos \phi) \be \cos \phi}{G'(\be \cos \phi)}
 -\frac{G(\be) \be}{G'(\be)} \right) \frac{G'(\be) G'(\be \cos \phi)}{G(\be) \be \sin^2 \phi}
\\
&=\frac{G'(\be) G'(\be \cos \phi)(H(\be) -H(\be \cos \phi))}{G(\be) \be \sin^2 \phi}.
\end{align*}
This is positive, since the assumption \eqref{aasfw} implies \eqref{aasGw}.
Hence it follows that
\begin{equation}\label{igjg}
I(\be,\phi) \int_0^{\phi} G''(\be \cos \tau) \cos \tau d\tau
 +J(\be,\phi) G'(\be\cos \phi) >0.
\end{equation}
Combining \eqref{sgni} and \eqref{igjg},
we obtain
\begin{equation*}
\be \frac{d}{d\be} \sqrt{\lam_k^o(\be)}>0.
\end{equation*}
Thus (ii) is verified, and the proof is complete.
\end{proof}

\section{Secondary bifurcations}\label{sec:sb}

In this section, we consider bifurcation points on ${\cS}^o$ and ${\cS}^e$.

\subsection{Nonexistence of bifurcation points on ${\cS}^e$}

The following lemma shows that no bifurcation point exists on ${\cS}^e$.

\begin{lem}\label{lem:ndse}
Assume \eqref{aasfw}. 
Then for every $(\lam,u) \in {\cS}^e$,
there is a neighborhood $\cN$ of $(\lam,u)$ in $(0,\infty) \times X_0$
such that ${\cS} \cap {\cN} ={\cS}^e \cap {\cN}$.
\end{lem}

\begin{proof}
By the assumption \eqref{aasfw} and Proposition~\ref{prop:pb},
we see that ${\cS}^e$ is the union of $C^1$ curves parametrized by $\lam$.
Therefore, according to (i) of Proposition~\ref{prop:bt},
we only need to show that $D(\lam_k^e (\be),\be,\be) \neq 0$ 
for $\be \in I \setminus \{0\}$.
Using \eqref{pqdb} and the fact that $\theta (\lam_k^e (\be),\be)=k\pi$ gives
\begin{align*}
D(\lam_k^e (\be),\be,\be) &=2P_\be (\lam_k^e (\be),\be) Q_\be (\lam_k^e (\be),\be)
\\
&=2\left\{ G'(\be) +a\frac{\be}{G'(\be)}
 \left(  \int_0^{k\pi} G'(\be \cos \tau) d\tau \right)
  \left( \int_0^{k\pi} G''(\be \cos \tau) \cos \tau d\tau \right) \right\}
\\
&\quad \times \frac{\be}{G'(\be)} \int_0^{k\pi} G''(\be \cos \tau) \cos \tau d\tau.
\end{align*}
Lemma~\ref{lem:gdpe} shows that this is positive,
and hence the lemma follows.
\end{proof}

\subsection{The number of bifurcation points on ${\cS}^o$}

We show that ${\cC}_{k,+}^o$ and ${\cC}_{k,-}^o$ have a unique bifurcation point,
provided that \eqref{aasfs} holds.
\begin{prop}\label{prop:bpco}
Assume \eqref{aasfs}.
Then for $k \in \bbn$, 
there exists a $C^1$ curve $\tilde {\cC}_{k,+}^o \subset {\cS}$ such that
$\tilde {\cC}_{k,+}^o$ intersects ${\cC}_{k,+}^o$ transversally 
at some point $(\lam^*_{k,+},u^*_{k,+})$.
Moreover, for each $(\lam,u) \in {\cC}_{k,+}^o$,
there is a neighborhood ${\cN}$ of $(\lam,u)$ such that
\begin{equation*}
{\cS} \cap \cN =\left\{
\begin{aligned}
&{\cC}_{k,+}^o \cap {\cN} && \mbox{if } (\lam,u) \neq (\lam^*_{k,+},u^*_{k,+}),
\\
&({\cC}_{k,+}^o \cup \tilde {\cC}_{k,+}^o) \cap {\cN}
 && \mbox{if } (\lam,u)=(\lam^*_{k,+},u^*_{k,+}).
\end{aligned}
\right.
\end{equation*}
The same assertion holds for ${\cC}_{k,-}^o$ in place of ${\cC}_{k,+}^o$.
\end{prop}

To prove this proposition, we examine the behavior of $D(\lam_k^o (\be),\be,-\be)$.
First we consider $P_\be (\lam_k^o (\be),\be)$.
\begin{lem}\label{lem:Pbnz}
If \eqref{aasfw} holds, 
then $(-1)^{k-1} P_\be (\lam_k^o (\be),\be)>0$ 
for all $k \in \bbn$ and $\be \in I \setminus \{0\}$.
\end{lem}

\begin{proof}
By definition, we have 
\begin{equation}\label{thph}
\theta(\lam_k^o (\be),\be)=\phi_k(\be).
\end{equation}
This together with \eqref{pqdb} and \eqref{ge0e} yields
\begin{align*}
P_\be(\lam_k^o (\be),\be) 
 &=\left. \frac{G'(\be \cos \phi) \be \cos \phi -G(\be \cos \phi)}{\be}
  \right|_{\phi=\phi_k(\be)}
\\
&\quad \left. +\left( \be \sin \phi
 +\frac{G(\be \cos \phi) \cos \phi}{G'(\be \cos \phi) \sin \phi} \right)
  \int_0^{\phi} G''(\be \cos \tau) \cos \tau d\tau \right|_{\phi=\phi_k(\be)}.
\end{align*}
It follows from \eqref{igjg0} that
\begin{align*}
&\left. \frac{G'(\be \cos \phi) \be}{G(\be \cos \phi)} 
 \int_0^{\phi}G'(\be \cos \tau)d\tau \right|_{\phi=\phi_k(\be)} \cdot P_\be(\lam_k^o (\be),\be)
\\
&=\left. \left( I(\be,\phi) \int_0^{\phi} G''(\be \cos \tau) \cos \tau d\tau
 +J(\be,\phi) G'(\be\cos \phi) \right) \right|_{\phi=\phi_k(\be)}.
\end{align*}
Therefore the desired inequality follows from \eqref{igjg}.
\end{proof}

Next we consider $Q_\be(\lam_k^o (\be),\be)$.
We note that \eqref{pqdb} and \eqref{thph} yield
\begin{equation*}
Q_\be(\lam_k^o (\be),\be)=R(\be,\phi_k(\be)),
\quad
R(\be,\phi):=-\sin \phi +\frac{\be \cos \phi}{G'(\be \cos \phi)}
 \int_0^\phi G''(\be \cos \tau) \cos \tau d\tau.
\end{equation*}

\begin{lem}\label{lem:qb0l}
For every $k \in \bbn$,
\begin{equation}\label{Qoas}
(-1)^{k-1} Q_\be(\lam_k^o (0),0)<0,
\qquad
\lim_{\be \to \pm \be_0} (-1)^{k-1} Q_\be(\lam_k^o (\be),\be) =\infty.
\end{equation}
\end{lem}

\begin{proof}
We see from \eqref{phpr1} that $Q_\be(\lam_k^o (0),0)=R(0,z_k)=-\sin z_k$.
Hence the first inequality of \eqref{Qoas} holds.

We examine the limit of $Q_\be(\lam_k^o (\be),\be)$ as $\be \to \pm \be_0$. 
By \eqref{Gdpas}, we have
\begin{align*}
\be \int_0^{\phi_k(\be)} G''(\be \cos \tau) \cos \tau d\tau
 &\ge \int_0^{\phi_k(\be)} |\be \cos \tau|
  \left\{ \frac{c}{(\be_0-|\be \cos \tau|)^{3/2}} -C \right\} d\tau
\\ 
&\ge c\int_{(k-1)\pi}^{\phi_k(\be)} 
  \frac{|\be \cos \tau|}{(\be_0-|\be \cos \tau|)^{3/2}} d\tau -C|\be| \phi_k(\be).
\end{align*}
Since
\begin{equation*}
\be_0-|\be \cos \tau| =\be_0 -|\be| +|\be| (1-|\cos \tau|)
 \le \be_0 -|\be| +\be_0 \sin^2 \tau,
\end{equation*}
we find that
\begin{align*}
\int_{(k-1)\pi}^{\phi_k(\be)} \frac{|\cos \tau|}{(\be_0-|\be \cos \tau|)^{3/2}}d\tau
 &\ge \int_{(k-1)\pi}^{\phi_k(\be)} \frac{|\cos \tau|}{(\be_0 -|\be| +\be_0 \sin^2 \tau)^{3/2}}d\tau
\\
&=\frac{1}{\sqrt{\be_0}(\be_0-|\be|)}
 \int_0^{\frac{\sqrt{\be_0}}{\sqrt{\be_0-|\be|}}|\sin \phi_k(\be)|} \frac{d\eta}{(1+\eta^2)^{3/2}},
\end{align*}
where we have used the change of variables 
$\eta=\sqrt{\be_0}|\sin \tau|/\sqrt{\be_0-|\be|}$.
Lemma~\ref{lem:phli} implies that
the integral on the right is bounded below by some positive constant.
Therefore
\begin{equation}\label{Qoes2}
\be \int_0^{\phi_k(\be)} G''(\be \cos \tau) \cos \tau d\tau
 \ge \frac{\tilde c}{\be_0-|\be|} -C|\be| \phi_k(\be),
\end{equation}
where $\tilde c>0$ is a constant.
From \eqref{Gpas}, we see that
\begin{equation}\label{Qoes1}
G'(\be \cos \phi_k(\be)) \le \frac{C}{\sqrt{\be_0-|\be \cos \phi_k(\be)|}}
 \le \frac{C}{\sqrt{\be_0-|\be|}}.
\end{equation}
Combining \eqref{Qoes2}, \eqref{Qoes1} and \eqref{phpr2} gives
\begin{align*}
&(-1)^{k-1} Q_\be(\lam_k^o (\be),\be) 
\\
&\ge (-1)^k \sin \phi_k(\be) +(-1)^{k-1} \cos \phi_k(\be)
 \left( \frac{\tilde c}{\sqrt{\be_0-|\be|}} -C|\be| \phi_k(\be) \sqrt{\be_0-|\be|}\right)
\\
&\to \infty \quad (\be \to \pm \be_0).
\end{align*}
We thus obtain \eqref{Qoas}.
\end{proof}

Lemma~\ref{lem:qb0l} shows that 
$Q_\be(\lam_k^o (\be^*),\be^*)=0$ for some $\be^* \in I \setminus\{0\}$.
In what follows, we fix such $\be^*$.
Set $\phi^*:=\phi_k(\be^*)$.
Since $Q_\be(\lam_k^o (\be),\be)=R(\be,\phi_k (\be))$, 
we see that the condition $Q_\be(\lam_k^o (\be^*),\be^*)=0$ is equivalent to
\begin{equation}\label{Qoez}
\frac{\be^* \cos \phi^*}{G'(\be^* \cos \phi^*)}
 \int_0^{\phi^*} G''(\be^* \cos \tau) \cos \tau d\tau =\sin \phi^*.
\end{equation}
We investigate the properties of $(\be^*,\phi^*)$.

\begin{lem}\label{lem:phes}
There holds 
\begin{equation*}
1-\frac{\be^* \sin \phi^*}{\cos \phi^*} \frac{d\phi_k}{d\be}(\be^*)>0.
\end{equation*}
\end{lem}

\begin{proof}
A direct computation yields
\begin{align*}
&\frac{1}{\be^*} I(\be^*,\phi^*) -\frac{\sin \phi^*}{\cos \phi^*}J(\be^*,\phi^*) \\
&=\frac{1}{\sin \phi^* \cos \phi^*} \int_0^{\phi^*} G'(\be^* \cos \tau)d\tau
 +G'(\be^* \cos \phi^*) +\frac{\be^* \sin \phi^*}{\cos \phi^*}
  \int_0^{\phi^*} G''(\be^* \cos \tau) \cos \tau d\tau
\\
&=\frac{1}{\sin \phi^* \cos \phi^*} \int_0^{\phi^*} G'(\be^* \cos \tau)d\tau
 +\frac{G'(\be^* \cos \phi^*)}{\cos^2 \phi^*}
\\
&>0,
\end{align*}
where we have used \eqref{Qoez} in deriving the second equality.
This together with \eqref{sgni} shows that
\begin{equation*}
1-\frac{\be^* \sin \phi^*}{\cos \phi^*} \frac{d\phi_k}{d\be}(\be^*)
=\frac{\be^*}{I(\be^*,\phi^*)}
 \left( \frac{1}{\be^*} I(\be^*,\phi^*) -\frac{\sin \phi^*}{\cos \phi^*}J(\be^*,\phi^*) \right) >0,
\end{equation*}
as claimed.
\end{proof}

We recall that the function $h(v)$ is given by \eqref{hdef}.
\begin{lem}\label{lem:bclv}
If \eqref{aasfs0} and \eqref{aasfw} hold,
then $h(\be^* \cos \phi^*)>0$.
\end{lem}
\begin{proof}
To obtain a contradiction, we suppose that $h(\be^* \cos \phi^*) \le 0$.
By Lemma~\ref{lem:gdpe},
we have
\begin{equation}\label{bclv1}
\be^* \int_0^{\phi^*} G''(\be^* \cos \tau) \cos \tau d\tau
 \ge \be^* \int_{(k-1)\pi}^{\phi^*} G''(\be^* \cos \tau) \cos \tau d\tau.
\end{equation}
According to (ii) of Lemma~\ref{lem:Gpro}, we can use \eqref{aashs0}.
Hence, from the assumption $h(\be^* \cos \phi^*) \le 0$,
we find that
\begin{equation*}
h(\be^* \cos \tau)<0 \quad \mbox{for all } \tau \in [(k-1)\pi,\phi^*).
\end{equation*}
This particularly implies that the function $G'(\be^* \cos \tau)/|\cos \tau|$ 
is decreasing on $[(k-1)\pi,\phi^*)$, since $(G'(v)/v)'=-G'(v)h(v)/v^2$.
Therefore
\begin{align*}
G''(\be^* \cos \tau) \be^* \cos \tau
 &=G'(\be^* \cos \tau)(1-h(\be^* \cos \tau))
\\
&>\frac{G'(\be^* \cos \phi^*)}{|\cos \phi^*|} |\cos \tau| 
 \quad \mbox{for all } \tau \in [(k-1)\pi,\phi^*).
\end{align*}
Plugging this into \eqref{bclv1},
we obtain
\begin{equation*}
\be^* \int_0^{\phi^*} G''(\be^* \cos \tau) \cos \tau d\tau
 >\frac{G'(\be^* \cos \phi^*)}{|\cos \phi^*|} 
  \int_{(k-1)\pi}^{\phi^*} |\cos \tau| d\tau
   =\frac{G'(\be^* \cos \phi^*) \sin \phi^*}{\cos \phi^*}.
\end{equation*}
This contradicts \eqref{Qoez}, and therefore $h(\be^* \cos \phi^*)>0$.
\end{proof}

\begin{lem}\label{lem:rbes1}
Let \eqref{aasfs0} and \eqref{aasfs1} hold and assume that $h(\be^*)>0$.
Then 
\begin{equation}\label{rbe1}
R_\be(\be^*,\phi^*) \be^* \sin \phi^* >0.
\end{equation}
\end{lem}

\begin{proof}
We note that the assumptions \eqref{aasfs0} and \eqref{aasfs1} 
give \eqref{aashs0} and \eqref{aashs1}.
By \eqref{aashs0} and the assumption $h(\be^*)>0$,
we have
\begin{equation}\label{bces2}
|\be^* \cos \tau | <v_0 \quad \mbox{for all } \tau \in [0,\phi^*].
\end{equation}
A direct computation yields
\begin{align*}
R_\be (\be,\phi)&=-\frac{G''(\be \cos \phi) \be \cos^2 \phi}{G'(\be \cos \phi)^2}
 \int_0^\phi G''(\be \cos \tau) \cos \tau d\tau
\\
&\quad +\frac{\cos \phi}{G'(\be \cos \phi)}
 \frac{\partial}{\partial \be} \left( \be \int_0^\phi G''(\be \cos \tau) \cos \tau d\tau \right).
\end{align*}
Since $(G''(v)v)' =-(G'(v)h(v))'+G''(v)$,
we have
\begin{align*}
\frac{\partial}{\partial \be} \left( G''(\be \cos \tau) \be \cos \tau \right)
 &=-\frac{d}{dv} (G'(v) h(v))|_{v=\be \cos \tau} \cdot \cos \tau +G''(\be \cos \tau) \cos \tau
\\
&=\frac{\D}{\D \tau} \left( G'(\be \cos \tau) h(\be \cos \tau)\right)
 \cdot \frac{\cos \tau}{\be \sin \tau} +G''(\be \cos \tau) \cos \tau,
\end{align*}
and hence
\begin{align}
R_\be (\be,\phi) &=\frac{h(\be \cos \phi)\cos \phi}{G'(\be \cos \phi)}
 \int_0^\phi G''(\be \cos \tau) \cos \tau d\tau
\nonumber
\\
&\quad +\frac{\cos \phi}{G'(\be \cos \phi) \be}
 \int_0^\phi \frac{\D}{\D \tau} \left( G'(\be \cos \tau) h(\be \cos \tau)\right)
  \cdot \frac{\cos \tau}{\sin \tau} d\tau.
\label{rbre}
\end{align}
Note that
\begin{align*}
&h(\be \cos \phi) \int_0^\phi G''(\be \cos \tau) \cos \tau d\tau
\\
&=(h(\be \cos \phi)-h(\be)) \int_0^\phi G''(\be \cos \tau) \cos \tau d\tau
 -\frac{h(\be)}{\be} \int_0^\phi \left( \frac{\D}{\D \tau} G'(\be \cos \tau)\right)
  \cdot \frac{\cos \tau}{\sin \tau} d\tau.
\end{align*}
Substituting this into \eqref{rbre}, we obtain
\begin{align*}
R_\be(\be,\phi) &=\frac{(h(\be \cos \phi)-h(\be)) \cos \phi}{G'(\be \cos \phi)} 
 \int_0^\phi G''(\be \cos \tau) \cos \tau d\tau 
\\
&\quad +\frac{\cos \phi}{G'(\be \cos \phi) \be}
 \int_0^\phi \frac{\D}{\D \tau} \big\{ G'(\be \cos \tau) (h(\be \cos \tau) -h(\be))
  \big\} \cdot \frac{\cos \tau}{\sin \tau} d\tau.
\end{align*}
We now apply integration by parts to the second term on the right.
Then, since 
\begin{equation*}
h(\be \cos \tau) -h(\be) =O(1-|\cos \tau|)=O(\sin^2 \tau)
\quad
\mbox{as } \tau \to l\pi, \ l \in \bbz,
\end{equation*}
we have
\begin{align*}
&\int_0^\phi \frac{\D}{\D \tau} \big\{ G'(\be \cos \tau) (h(\be \cos \tau) -h(\be))
 \big\} \cdot \frac{\cos \tau}{\sin \tau} d\tau
\\
&=\frac{G'(\be \cos \phi) (h(\be \cos \phi) -h(\be)) \cos \phi}{\sin \phi}
 +\int_0^\phi \frac{G'(\be \cos \tau)(h(\be \cos \tau) -h(\be))}{\sin^2 \tau} d\tau.
\end{align*}
This together with \eqref{Qoez} shows that
\begin{align*}
&R_\be (\be^*,\phi^*) \be^* \sin \phi^*
\\
&=\left( \frac{\be^* \sin \phi^* \cos \phi^*}{G'(\be^* \cos \phi^*)} 
 \int_0^{\phi^*} G''(\be^* \cos \tau) \cos \tau d\tau
  +\cos^2 \phi^* \right) (h(\be^* \cos \phi^*)-h(\be^*))
\\
&\quad +\frac{\sin \phi^* \cos \phi^*}{G'(\be^* \cos \phi^*)}
 \int_0^{\phi^*} \frac{G'(\be^* \cos \tau)(h(\be^* \cos \tau) -h(\be^*))}{\sin^2 \tau} d\tau
\\
&=h(\be^* \cos \phi^*)-h(\be^*)
 +\frac{\sin \phi^* \cos \phi^*}{G'(\be^* \cos \phi^*)}
  \int_0^{\phi^*} \frac{G'(\be^* \cos \tau)(h(\be^* \cos \tau) -h(\be^*))}{\sin^2 \tau} d\tau.
\end{align*}
We see from \eqref{aashs1} and \eqref{bces2} that the right-hand side is positive,
and \eqref{rbe1} is proved.
\end{proof}

\begin{lem}\label{lem:rbes2}
Let \eqref{aasfs0}, \eqref{aasfs2} and \eqref{aasfw} hold
and assume that $h(\be^*) \le 0$.
Then
\begin{equation}\label{rbe2}
R_\be (\be^*,\phi^*) \be^* \sin \phi^* >h(\be^* \cos \phi^*).
\end{equation}
\end{lem}

\begin{proof}
We see from (ii) of Lemma~\ref{lem:Gpro} that
\eqref{aashs0} and \eqref{aashs2} are satisfied.
By \eqref{aashs0}, Lemma~\ref{lem:bclv} and the assumption $h(\be^*) \le 0$, 
we can take $\tau_0 \in [0,\phi^*-(k-1)\pi)$ such that $|\be^*| \cos \tau_0=v_0$.
Put
\begin{equation*}
{\cI}:=((k-1)\pi+\tau_0,\phi^*] \cup \bigcup_{j=1}^{k-1} ((j-1)\pi+\tau_0, j\pi -\tau_0).
\end{equation*}
Then, from \eqref{aashs0}, we have 
\begin{equation}\label{bces}
|\be^* \cos \tau| \left\{
\begin{aligned}
&<v_0 && \mbox{if } \tau \in {\cI},
\\
&\ge v_0 && \mbox{if } \tau \in [0,\phi^*] \setminus {\cI},
\end{aligned}
\right.
\end{equation}
and
\begin{equation}\label{shbc}
h(\be^*\cos \tau) \left\{
\begin{aligned}
&>0 && \mbox{if } \tau \in {\cI},
\\
&=0 && \mbox{if } \tau \in \partial {\cI} \setminus \{ \phi^*\}.
\end{aligned}
\right.
\end{equation}
We estimate the second term on the right of \eqref{rbre}.
Using integration by parts and \eqref{shbc},
we have
\begin{align*}
&\int_{\cI} \frac{\D}{\D \tau} \left( G'(\be^* \cos \tau) h(\be^* \cos \tau)\right)
  \cdot \frac{\cos \tau}{\sin \tau} d\tau
\\
&=\frac{G'(\be^* \cos \phi^*) h(\be^* \cos \phi^*)\cos \phi^*}{\sin \phi^*} 
 +\int_{\cI} \frac{G'(\be^* \cos \tau) h(\be^* \cos \tau)}{\sin^2 \tau} d\tau
\\
&>\frac{G'(\be^* \cos \phi^*) h(\be^* \cos \phi^*)\cos \phi^*}{\sin \phi^*}.
\end{align*}
Furthermore, \eqref{aashs2} and \eqref{bces} yield
\begin{align*}
&\int_{[0,\phi^*] \setminus {\cI}}
 \frac{\D}{\D \tau} \left( G'(\be^* \cos \tau) h(\be^* \cos \tau)\right)
  \cdot \frac{\cos \tau}{\sin \tau} d\tau
\\
&=-\int_{[0,\phi^*] \setminus {\cI}}
 \left. \left\{ v \frac{d}{dv} (G'(v) h(v)) \right\} \right|_{v=\be^* \cos \tau} d\tau
\\
&\ge 0.
\end{align*}
Plugging these inequalities into \eqref{rbre} and using \eqref{Qoez},
we obtain
\begin{align*}
&R_\be (\be^*,\phi^*) \be^* \sin \phi^*
\\
&>\frac{h(\be^* \cos \phi^*) \be^* \sin \phi^* \cos \phi^*}{G'(\be^* \cos \phi^*)}
 \int_0^{\phi^*} G''(\be^* \cos \tau) \cos \tau d\tau
  +h(\be^* \cos \phi^*) \cos^2 \phi^*
\\
&=h(\be^* \cos \phi^*),
\end{align*}
which proves the lemma.
\end{proof}

\begin{lem}\label{ddbp}
Assume \eqref{aasfs}.
Then 
\begin{equation}\label{dQnz}
(-1)^{k-1} \sgn (\be^*) \left. \frac{d}{d\be}Q_\be(\lam_k^o (\be),\be) \right|_{\be=\be^*}>0.
\end{equation}
\end{lem}

\begin{proof}
Since $Q_\be(\lam_k^o (\be),\be)=R(\be,\phi_k(\be))$, we have
\begin{equation*}
\left. \frac{d}{d\be}Q_\be(\lam_k^o (\be),\be) \right|_{\be=\be^*}
  =R_\be (\be^*,\phi^*) +R_\phi (\be^*,\phi^*) \frac{d\phi_k}{d\be}(\be^*).
\end{equation*}
To estimate this, we compute $R_\phi (\be^*,\phi^*)$.
A direct computation gives
\begin{align*}
R_\phi (\be,\phi)&=-\cos \phi
 -\frac{(G'(\be \cos \phi)-G''(\be \cos \phi) \be \cos \phi) \be \sin \phi}{G'(\be \cos \phi)^2}
  \int_0^\phi G''(\be \cos \tau) \cos \tau d\tau
\\  
&\quad +\frac{\be \cos \phi}{G'(\be \cos \phi)} \cdot G''(\be \cos \phi) \cos \phi
\\
&=-\left( \cos \phi +\frac{\be \sin \phi}{G'(\be \cos \phi)}
 \int_0^\phi G''(\be \cos \tau) \cos \tau d\tau \right) h(\be \cos \phi).
\end{align*}
From \eqref{Qoez}, we find that
\begin{equation}\label{rphe}
R_\phi (\be^*,\phi^*)
  =-\frac{h(\be^* \cos \phi^*)}{\cos \phi^*}.
\end{equation}

We consider the two cases: $h(\be^*)>0$ and $h(\be^*) \le 0$.
We first consider the latter case.
Lemma~\ref{lem:rbaf} shows that \eqref{aasfw} is satisfied,
and hence we can apply Lemmas~\ref{lem:bclv} and \ref{lem:rbes2}.
From \eqref{rbe2} and \eqref{rphe}, 
we have
\begin{align*}
\be^* \sin \phi^* \cdot \left. \frac{d}{d\be}Q_\be(\lam_k^o (\be),\be) \right|_{\be=\be^*}
 >h(\be^* \cos \phi^*)
  \left( 1 -\frac{\be^* \sin \phi^*}{\cos \phi^*} \frac{d\phi_k}{d\be}(\be^*) \right).
\end{align*}
Lemmas~\ref{lem:phes} and \ref{lem:bclv} show that the right-hand side is positive.
Since the sign of $\be^* \sin \phi^*$ coincides with that of $(-1)^{k-1} \sgn (\be^*)$,
we obtain \eqref{dQnz}.

Let us consider the other case $h(\be^*)>0$. 
Then we have $(-1)^{k-1} \sgn (\be^*) R_\be (\be^*,\phi^*)>0$ by Lemma~\ref{lem:rbes1}.
Moreover, we see from Lemma~\ref{lem:bclv} and \eqref{rphe} that 
$(-1)^{k-1}R_\phi (\be^*,\phi^*)<0$.
Therefore it is sufficient to show that
\begin{equation*}
\sgn (\be^*) \frac{d\phi_k}{d\be}(\be^*) \le 0.
\end{equation*}
To prove this, we define
\begin{equation*}
\tilde h(v):=\left\{
\begin{aligned}
&\frac{G'(v)v}{G(v)} && (v \neq 0),
\\
&1 && (v=0).
\end{aligned}
\right.
\end{equation*}
By (i) of Lemma~\ref{lem:Gpro}, 
we see that $\tilde h \in C(I) \cap C^2 (I \setminus \{0\})$ and $G\tilde h \in C^1(I)$.
Notice that
\begin{equation*}
\frac{d}{dv} \left( \frac{G(v)^2 \tilde h'(v)}{G'(v)} \right) =-G(v)h'(v)
 \quad \mbox{for } v \in I \setminus \{0\}.
\end{equation*}
\eqref{aashs1} shows that the right-hand side of this equality is nonnegative if $v \in (v_0,v_0)$,
and hence
\begin{equation*}
\sgn(v) \frac{G(v)^2 \tilde h'(v)}{G'(v)}
 \ge \lim_{v \to 0} \sgn(v) \frac{G(v)^2 \tilde h'(v)}{G'(v)}
   =0 \quad \mbox{for } v \in (v_0,v_0) \setminus \{0\}.
\end{equation*}
From this we have
\begin{equation}\label{edth}
\sgn(v) \tilde h(v) \mbox{ is nondecreasing in } (v_0,v_0) \setminus \{0\}.
\end{equation}
Since $G'(v)+G''(v)v=(G(v)\tilde h(v))'$, we see that $J(\be,\phi)$ is written as
\begin{align*}
J(\be,\phi)
&=\tilde h(\be \cos \phi) \int_0^\phi G'(\be \cos \tau)d\tau 
 +\int_0^\phi \frac{\D}{\D \tau}
  \big( G(\be \cos \tau) \tilde h(\be \cos \tau) \big) \cdot \frac{1}{\be \sin \tau} d\tau
\\
&=(\tilde h(\be \cos \phi) -\tilde h(\be)) \int_0^\phi G'(\be \cos \tau)d\tau 
\\
&\quad +\int_0^\phi \frac{\D}{\D \tau} \big\{ G(\be \cos \tau) (\tilde h(\be \cos \tau)
 -\tilde h(\be)) \big\} \cdot \frac{1}{\be \sin \tau} d\tau.
\end{align*}
Integrating by parts shows that
the second term on the right is computed as
\begin{align*}
&\int_0^\phi \frac{\D}{\D \tau}
 \big\{ G(\be \cos \tau) (\tilde h(\be \cos \tau) -\tilde h(\be)) \big\} \cdot \frac{1}{\be \sin \tau} d\tau
\\
&=(\tilde h(\be \cos \phi) -\tilde h(\be)) \frac{G(\be \cos \phi)}{\be \sin \phi} 
 +\int_0^\phi (\tilde h(\be \cos \tau) -\tilde h(\be))
  \frac{G(\be \cos \tau) \cos \tau}{\be \sin^2 \tau} d\tau.
\end{align*}
Therefore
\begin{align}
J(\be^*,\phi^*) &=(\tilde h(\be^* \cos \phi^*) -\tilde h(\be^*))
 \left( \int_0^{\phi^*} G'(\be^* \cos \tau)d\tau
  +\frac{G(\be^* \cos \phi^*)}{\be^* \sin \phi^*}\right)
\nonumber
\\
&\quad +\int_0^{\phi^*} (\tilde h(\be^* \cos \tau) -\tilde h(\be^*))
 \frac{G(\be^* \cos \tau) \cos \tau}{\be^* \sin^2 \tau} d\tau.
\label{Jbpa}
\end{align}
By \eqref{aashs0} and the assumption $h(\be^*)>0$,
we have $|\be^* \cos \tau|<v_0$ for all $\tau \in [0,\phi^*]$.
It follows from \eqref{edth} that the right-hand side of \eqref{Jbpa} is nonpositive,
and hence $J(\be^*,\phi^*) \le 0$.
This together with \eqref{sgni} shows that
\begin{equation*}
\sgn (\be^*) \frac{d\phi_k}{d\be}(\be^*) 
 =\frac{\sgn (\be^*)}{I(\be^*,\phi^*)} \cdot J(\be^*,\phi^*) \le 0.
\end{equation*}
This completes the proof.
\end{proof}

We can now prove Proposition~\ref{prop:bpco}.
\begin{proof}[Proof of Proposition~\ref{prop:bpco}]
By \eqref{PQsy}, we have
\begin{equation*}
D(\lam_k^o (\be),\be,-\be)=2P_\be (\lam_k^o (\be),\be)Q_\be (\lam_k^o (\be),\be).
\end{equation*}
In particular, Lemma~\ref{lem:Pbnz} yields
\begin{equation*}
{\cD}:=\{ \be \in I \setminus \{0\}; D(\lam_k^o (\be),\be,-\be)=0\}
 =\{ \be \in I \setminus \{0\}; Q_\be (\lam_k^o (\be),\be)=0\}.
\end{equation*}
We see from Lemmas~\ref{lem:Pbnz} and \ref{lem:qb0l} that
\begin{equation*}
D(\lam_k^o (\be),\be,-\be) \left\{
\begin{aligned}
&<0 && \mbox{if } |\be| \mbox{ is small},
\\
&>0 && \mbox{if } |\be| \mbox{ is close to } \be_0.
\end{aligned}
\right.
\end{equation*}
Furthermore, Lemmas~\ref{lem:Pbnz} and \ref{ddbp} show that
for any $\be^* \in {\cD}$,
\begin{equation}\label{dDba}
\sgn (\be^*) \left. \frac{d}{d\be} D(\lam_k^o (\be),\be,-\be) \right|_{\be=\be^*} 
 =2\sgn (\be^*) P_\be (\lam_k^o (\be^*),\be^*)
  \cdot \left. \frac{d}{d\be} Q_\be (\lam_k^o (\be),\be) \right|_{\be=\be^*} >0.
\end{equation}
Therefore there exist $\be_+^* \in (0,\be_0)$ and $\be_-^* \in (-\be_0,0)$ such that
\begin{equation}\label{cDtp}
{\cD}=\{ \be_+^*,\be_-^*\}.
\end{equation}

From Lemma~\ref{lem:rbaf}, 
we know that \eqref{aasfw} is satisfied under the assumption \eqref{aasfs}.
Hence (ii) of Proposition~\ref{prop:pb} shows that 
${\cC}_{k,+}^o$ and ${\cC}_{k,-}^o$ are written as
\begin{equation}\label{copl}
{\cC}_{k,+}^o=\{ (\lam,u^o_k(\cdot,\be_+^o(\lam)))\}_{\lam \in (\lam_{2k-1},\infty)},
\qquad
{\cC}_{k,-}^o=\{ (\lam,u^o_k(\cdot,\be_-^o(\lam)))\}_{\lam \in (\lam_{2k-1},\infty)},
\end{equation}
where $\be_+^o(\lam)$ (resp. $\be_-^o(\lam)$) 
is the inverse of the function $(0,\be_0) \ni \be \mapsto \lam_k^o(\be)$
(resp. $(-\be_0,0) \ni \be \mapsto \lam_k^o(\be)$).
Then, by \eqref{dDba}, \eqref{cDtp} and \eqref{sdlk}, we have
\begin{equation}\label{dDab}
\left\{
\begin{gathered}
D(\lam,\be_\pm^o(\lam),-\be_\pm^o(\lam))=0
 \mbox{ if and only if } \lam=\lam_k^o(\be_\pm^*),
\\
\left. \frac{d}{d\lam} D(\lam,\be_\pm^o(\lam),-\be_\pm^o(\lam))
 \right|_{\lam=\lam_k^o(\be_\pm^*)}
  =\left( \frac{d\lam_k^o}{d\be}(\be_\pm^*) \right)^{-1} \cdot 
   \left. \frac{d}{d\be} D(\lam_k^o (\be),\be,-\be) \right|_{\be=\be_\pm^*}>0.
\end{gathered}
\right.
\end{equation}
Thus we obtain the desired conclusion by applying (i) and (ii) of Proposition~\ref{prop:bt}.
\end{proof}

\begin{remk}\label{rem:sybp}
Since \eqref{PQsy} and \eqref{lksy} give
$Q_\be (\lam_k^o (-\be),-\be)=-Q_\be (\lam_k^o (\be),\be)$,
we infer that $\be_-^*=-\be_+^*$.
Hence the bifurcation points 
$(\lam_k^o(\be_+^*),u^o_k(\cdot,\be_+^*)) \in {\cC}_{k,+}^o$
and $(\lam_k^o(\be_-^*),u^o_k(\cdot,\be_-^*)) \in {\cC}_{k,-}^o$ 
obtained in Proposition~\ref{prop:bpco} satisfy
\begin{equation*}
(\lam_k^o(\be_-^*),u^o_k(\cdot,\be_-^*))
 =(\lam_k^o(\be_+^*),-u^o_k(\cdot,\be_+^*)).
\end{equation*}
\end{remk}

\subsection{Remark on the assumption \eqref{aasfs}}

At the end of this section, 
we observe that Proposition~\ref{prop:bpco} is still true for $k=1$ 
if we drop the assumption \eqref{aasfs1}.
\begin{prop}
Under the assumptions \eqref{aasfs0}, \eqref{aasfs2} and \eqref{aasfw},
Proposition~\ref{prop:bpco} holds for $k=1$.
\end{prop}

\begin{proof}
Let $\be^* \in I \setminus \{0\}$ satisfy $Q(\lam_1^o (\be^*),\be^*)=0$.
If \eqref{dQnz} is satisfied for $k=1$,
then the proposition can be proved in the same way as Proposition~\ref{prop:bpco}.
Therefore we only need to show that
\begin{equation}\label{hban}
h(\be^*) \le 0,
\end{equation}
since this enables us to apply Lemma~\ref{lem:rbes2} to obtain \eqref{dQnz}.
Contrary to \eqref{hban}, suppose that $h(\be^*)>0$.
Then \eqref{aashs0} yields 
\begin{equation*}
h(\be^* \cos \tau)>0 \quad \mbox{for all } \tau \in [0,\phi^*],
\end{equation*}
where $\phi^*:=\phi_1(\be^*) \in (0,\pi/2)$.
From this and the fact that $(G'(v)/v)'=-G'(v)h(v)/v^2$,
we see that the function $G'(\be^* \cos \tau)/\cos \tau$ 
is increasing on $[0,\phi^*]$.
Hence
\begin{align*}
\be^* \int_0^{\phi^*} G''(\be^* \cos \tau) \cos \tau d\tau
 &=\int_0^{\phi^*} G'(\be^* \cos \tau) (1-h(\be^* \cos \tau)) d\tau
\\
&<\int_0^{\phi^*} \frac{G'(\be^* \cos \phi^*)}{\cos \phi^*} \cos \tau d\tau
\\
&=\frac{G'(\be^* \cos \phi^*) \sin \phi^*}{\cos \phi^*}.
\end{align*}
This gives $Q(\lam_1^o (\be^*),\be^*)<0$, a contradiction.
Therefore we obtain \eqref{hban}, and the proof is complete.
\end{proof}

\section{Proof of Theorem~\ref{mthm}}\label{sec:pt}


To prove Theorem~\ref{mthm},
we compute the Morse index of solutions on ${\cS}^e$ and ${\cS}^o$.
We write $\lam^*$ for the number $\lam_{k,+}^*$ obtained in Proposition~\ref{prop:bpco}.
\begin{prop}\label{prop:mi}
For $k \in \bbn$, 
the following hold.
\begin{enumerate}[(i)]
\item
Let \eqref{aasfw} hold and let $(\lam,u) \in {\cC}_k^e$. 
Then $u$ is nondegenerate and $i(u)=2k$.
\item
Let \eqref{aasfs} hold and let $(\lam,u) \in {\cC}_k^o$. 
Then $u$ is nondegenerate unless $\lam \neq \lam^*$ and 
\begin{equation*}
i(u)=\left\{
\begin{aligned}
&2k-1 && (\lam<\lam^*),
\\
&2k-2 && (\lam \ge \lam^*).
\end{aligned}
\right.
\end{equation*}
\end{enumerate}
\end{prop}

In what follows, we fix $k \in \bbn$.
For $n \in \bbn \cup \{0\}$,
let $\mu_n^o(\be)$ (resp. $\mu_n^e(\be)$) 
denote the $(n+1)$-th largest eigenvalue of \eqref{llevp} 
for $(\lam,u)=(\lam_k^o(\be),u^o_k(\cdot,\be)) \in {\cC}_k^o$
(resp. $(\lam,u)=(\lam_k^e(\be),u^e_k(\cdot,\be)) \in {\cC}_k^e$).
We see from Lemma~\ref{lem:ipmi} that 
$\mu_n^o(\be)$ and $\mu_n^e(\be)$ are continuous with respect to $\be$.
In the following two lemmas, 
we give basic estimates of $\mu_n^o(\be)$ and $\mu_n^e(\be)$.
\begin{lem}\label{lem:evaz}
There hold $\mu_{2k-2}^o(0)>0$ and $\mu_{2k-1}^e(0)>0$.
\end{lem}

\begin{proof}
It is elementary to show that 
the $(n+1)$-th eigenvalue $\mu_n$ of \eqref{llevp} for $u=0$
is given by 
\begin{equation*}
\mu_{2k-2}=\lam f'(0)-\{(k-1)\pi \}^2,
\qquad
\mu_{2k-1}=\lam f'(0)-z_k^2
\quad
(k \in \bbn).
\end{equation*}
Hence
\begin{equation*}
\mu_{2k-2}^o(0)>\mu_{2k-1}^o(0)=\lam_{2k-1} f'(0)-z_k^2=0,
\qquad
\mu_{2k-1}^e(0)>\mu_{2k}^e(0)=\lam_{2k} f'(0)-(k\pi)^2=0,
\end{equation*}
as desired.
\end{proof}

\begin{lem}\label{lem:eves}
Assume that \eqref{aasfw} holds.
Then $\mu_{2k-1}^o(\be)<0$ and $\mu_{2k}^e(\be)<0$
for all $\be \in I \setminus \{0\}$.
\end{lem}

\begin{proof}
Let $Z(w)$ denote the number of zeros of a function $w$ in $(-1,1) \setminus \{0\}$.
By Lemma~\ref{lem:Mies},
it suffices to show that
\begin{gather*}
Z(u^o_k(\cdot,\be))=2k-2,
\quad
Z(u^e_k(\cdot,\be))=2k,
\\
u^o_k(-0,\be)u^o_k(+0,\be)<0,
\quad 
u^e_k(-0,\be)u^e_k(+0,\be)>0
\end{gather*}
for $\be \in I \setminus \{0\}$.
To derive these,
we recall that any $u \in {\cS}_\lam$ is written as
\begin{equation*}
u(x)=\left\{
\begin{aligned}
&U\left( \sqrt{\lambda} (x+1),\be_1 \right)
 =G\left( \be_1 \cos \Theta \left( \sqrt{\lambda} (x+1),\be_1\right) \right)
  && \mbox{for } x \in [-1,0), \\
&U\left( \sqrt{\lambda} (1-x),\be_2 \right)
 =G\left( \be_2 \cos \Theta \left( \sqrt{\lambda} (1-x),\be_2\right) \right)
  && \mbox{for } x \in (0,1],
\end{aligned}
\right.
\end{equation*}
where $\be_1=G(u(-1))$ and $\be_2=G(u(1))$.
This implies that if $\be_1 \be_2 \neq 0$ and
\begin{equation*}
\left( m_j-\frac{1}{2}\right) \pi< \theta (\lam,\be_j)
 =\Theta (\sqrt{\lam},\be_j)
  \le \left( m_j+\frac{1}{2}\right) \pi \quad \mbox{for some } m_j \in \bbn \cup \{0\}, 
\end{equation*}
then $Z(u)=m_1+m_2$ and $\sgn(u(-0)u(+0))=\sgn ((-1)^{m_1+m_2} \be_1 \be_2)$.
Since we know that $\theta(\lam_k^o (\be),\be)=\phi_k(\be) \in ((k-1)\pi,(k-1/2)\pi)$ 
and $\theta(\lam_k^e (\be),\be)=k\pi$,
we have
\begin{gather*}
Z(u^o_k(\cdot,\be))=(k-1)+(k-1)=2k-2,
\quad
Z(u^e_k(\cdot,\be))=k+k=2k,
\\
\sgn(u^o_k(-0,\be)u^o_k(+0,\be))=\sgn \left( (-1)^{(k-1)+(k-1)} \cdot (-\be^2)\right)
 =\sgn(-\be^2)<0,
\\
\sgn(u^e_k(-0,\be)u^e_k(+0,\be))=\sgn \left( (-1)^{k+k} \cdot \be^2\right)
 =\sgn( \be^2)>0
\end{gather*}
for $\be \in I \setminus \{0\}$.
Therefore the lemma follows.
\end{proof}

Let us show Proposition~\ref{prop:mi}.
\begin{proof}[Proof of Proposition~\ref{prop:mi}]
First we prove (i).
Lemma \ref{lem:evaz} shows that $\mu_{2k-1}^e(\be)$ is positive if $|\be|$ is small enough.
As shown in the proof of Lemma~\ref{lem:ndse},
we know that $D(\lam_k^e(\be),\be,\be) \neq 0$ for $\be \in I \setminus \{0\}$.
From this and Lemma~\ref{lem:ndcD}, we see that $\mu_{2k-1}^e(\be)$ never vanishes.
Therefore $\mu_{2k-1}^e(\be)>0$ for all $\be \in I$.
Thus (i) is verified by combining this with Lemma~\ref{lem:eves}.

Next we prove (ii).
We recall that \eqref{cDtp} holds.
Hence Lemma~\ref{lem:ndcD} gives
\begin{gather}
\mu_n^o(\be) \neq 0 \mbox{ for all } n \in \bbn \cup \{0\}
 \mbox{ and } \be \in I \setminus \{0,\be_+^*,\be_-^*\},
\label{ednv}
\\
\mu_{n_+}^o(\be_+^*)=\mu_{n_-}^o(\be_-^*)=0
 \mbox{ for some } n_+,n_- \in \bbn \cup \{0\}.
\nonumber
\end{gather}
Moreover, Lemma~\ref{lem:evaz} shows that
\begin{equation*}
\mu_n^o(\be) \mbox{ is positive if } |\be| \mbox{ is small and } n \le 2k-2.
\end{equation*}
Combining these with Lemma~\ref{lem:eves},
we deduce that
\begin{gather}
\mu_{2k-2}^o(\be_+^*)=\mu_{2k-2}^o(\be_-^*)=0,
\label{cevv}
\\
\mu_{2k-3}^o(\be) >0 \quad \mbox{for all } \be \in I,
 \mbox{ provided } k \ge 2.
\label{revp}
\end{gather}
To investigate the behavior of $\mu_{2k-2}^o(\be)$,
we apply (iii) of Proposition~\ref{prop:bt}.
For this purpose,
we use the parametrization of ${\cC}_{k,+}^o$ and ${\cC}_{k,+}^o$ as in \eqref{copl}.
Then $\mu_n^o(\be_\pm^o (\lam))=\mu_n(u^o_k(\cdot,\be_\pm^o (\lam))$.
By \eqref{cevv}, we can apply \eqref{evdf} 
for $n=2k-2$, $u(\cdot,\lam)=u^o_k(\cdot,\be_\pm^o (\lam))$,
$\be_1(\lam)=\be_\pm^o (\lam)$ and $\be_2(\lam)=-\be_\pm^o (\lam)$ to obtain
\begin{equation*}
\sgn \left( \left. \frac{d}{d\lam} \mu_{2k-2}^o(\be_\pm^o (\lam))
 \right|_{\lam=\lam_k^o(\be_\pm^*)} \right) =-\sgn \left( \left. \frac{d}{d\lam} 
  D(\lam,\be_{\pm}^o(\lam),-\be_{\pm}^o(\lam)) \right|_{\lam=\lam_k^o(\be_\pm^*)} \right).
\end{equation*}
According to \eqref{dDab}, we know that the right-hand side is negative.
Therefore it follows from \eqref{ednv} that
\begin{equation*}
\mu_{2k-2}^o(\be_\pm^o (\lam)) \left\{
\begin{aligned}
&>0 && \mbox{if } \lam<\lam_k^o(\be_\pm^*),
\\
&<0 && \mbox{if } \lam>\lam_k^o(\be_\pm^*).
\end{aligned}
\right.
\end{equation*}
Combining this with Lemma~\ref{lem:eves} and \eqref{revp},
we conclude that
\begin{equation*}
i(u^o_k(\cdot,\be_\pm^o (\lam))= \left\{
\begin{aligned}
&2k-1 && \mbox{if } \lam<\lam_k^o(\be_\pm^*),
\\
&2k-2 && \mbox{if } \lam \ge \lam_k^o(\be_\pm^*).
\end{aligned}
\right.
\end{equation*}
As noted in Remark~\ref{rem:sybp},
we know that $\lam_k^o(\be_+^*)=\lam_k^o(\be_-^*)$.
Consequently, (ii) is proved.
\end{proof}

We are now in a position to prove Theorem~\ref{mthm}.
\begin{proof}[Proof of Theorem~\ref{mthm}]
We note that by the assumption \eqref{aasfs}, the condition \eqref{aasfw} is satisfied.
We put ${\cC}_{2k-1}={\cC}_{k,+}^o$,  ${\cC}_{2k}={\cC}_{k,+}^e$.
Then (i) and (ii) follow immediately from 
Proposition~\ref{prop:pb} and Lemma~\ref{lem:ndse}.
Moreover, (iii) and (iv) are direct consequences of 
Propositions~\ref{prop:bpco} and \ref{prop:mi} and Remark~\ref{rem:sybp}.
Therefore the proof is complete.
\end{proof}


\section*{Acknowledgements}

The author would like to thank Professors Satoshi Tanaka and Masahito Ohta
for calling his attention to the references 
\cite{AN13,AG18,S90}.
This work is supported in part by the Grant-in-Aid for Early-Career Scientists 19K14574,
Japan Society for the Promotion of Science.


\appendix

\section{Relation between \eqref{aasfs} and \eqref{aasfw}}\label{appendixA}

In this section,
we show that \eqref{aasfs} implies \eqref{aasfw}.
\begin{lem}\label{lem:rbaf}
If \eqref{aasfs} holds, then \eqref{aasfw} is satisfied.
\end{lem}

\begin{proof}
Assume \eqref{aasfs}.
First we show that
\begin{equation}\label{aasfsw}
\frac{f'(u)F(u)}{f(u)^2}<1 \quad \mbox{for } u \in (-1,1) \setminus \{0\}.
\end{equation}
We see from \eqref{aasfs1} that for $u \in (-u_0,u_0) \setminus \{0\}$,
\begin{equation*}
\frac{f'(u)F(u)}{f(u)^2} <\lim_{s \to 0} \frac{f'(s)F(s)}{f(s)^2}
 =f'(0) \lim_{s \to 0} \frac{2f(s)}{2f(s)f'(s)}=1,
\end{equation*}
where we have used L'Hopital's rule in deriving the first equality.
Moreover, \eqref{aasfs0} yields $f'(u)F(u)/f(u)^2 \le 0$ 
for $u \in (-1,-u_0] \cup [u_0,1)$.
Therefore \eqref{aasfsw} holds.

Now we derive \eqref{aasfw}.
Notice that
\begin{equation*}
1-\frac{f'(u)F(u)}{f(u)^2} =\frac{d}{du} \left( \frac{F(u)}{f(u)} -u \right).
\end{equation*}
From this and the fact that $F(u)/f(u) \to 0$ ($u \to 0$), we have
\begin{equation*}
\int_0^u \left( 1-\frac{f'(s)F(s)}{f(s)^2}\right) ds =\frac{F(u)}{f(u)} -u
 =\frac{F(u)}{f(u)} \left( 1-\frac{f'(u)u}{f(u)}\right) -u\left( 1-\frac{f'(u)F(u)}{f(u)^2}\right).
\end{equation*}
Hence
\begin{equation*}
1-\frac{f'(u)u}{f(u)} =\frac{f(u)u}{F(u)} \left( 1-\frac{f'(u)F(u)}{f(u)^2}\right)
 +\frac{f(u)}{F(u)} \int_0^u \left( 1-\frac{f'(s)F(s)}{f(s)^2}\right) ds.
\end{equation*}
\eqref{aasfw} then follows immediately from \eqref{aasfsw}.
\end{proof}

\section{Proofs of Lemmas~\ref{lem:Gpro} and \ref{lem:bt1}}\label{appendixB}

This section provides the proofs of Lemmas~\ref{lem:Gpro} and \ref{lem:bt1}.
\begin{proof}[Proof of Lemma~\ref{lem:Gpro}]
By \eqref{basf}, we see that the function given by \eqref{Ginv} 
belongs to $C^3((-1,1) \setminus \{0\})$ 
and its derivative does not vanish in $(-1,1) \setminus \{0\}$.
Hence the inverse function theorem shows that $G \in C^3(I \setminus \{0\})$.

Let $u=G(v)$.
Differentiating the equality $F(u)=v^2$,
we have
\begin{gather}
G'(v)=\frac{v}{f(u)}=\frac{\sgn (u)\sqrt{F(u)}}{f(u)},
\qquad
G''(v)=\frac{1}{f(u)} \left( 1-\frac{f'(u)F(u)}{f(u)^2}\right),
\label{dGfo}
\\
\frac{d}{dv} (G'(v)h(v))=-G'''(v)v=\frac{\sqrt{F(u)}}{f(u)}
 \frac{d}{du} \left( \frac{f'(u)F(u)^{\frac{3}{2}}}{f(u)^3}\right),
\label{dG3fo}
\end{gather}
provided that $v \neq 0$.
Since Taylor's theorem yields
\begin{gather*}
f'(s)=f'(0)+f''(0)s+o(s),
\qquad
f(s)=f'(0)s+\frac{1}{2}f''(0)s^2+o(s^2),
\\
F(s)=f'(0)s^2 +\frac{1}{3}f''(0)s^3+o(s^3)
\quad \mbox{as } s \to 0,
\end{gather*}
we deduce that
\begin{equation}\label{Gdao}
\lim_{v \to 0}G'(v)=\frac{1}{\sqrt{f'(0)}},
\qquad
\lim_{v \to 0}G''(v)=-\frac{f''(0)}{3f'(0)^2}.
\end{equation}
Therefore $G \in C^2(I)$.
Note that $G''(v)v$ is written as $G''(v)v=G'(v)(1-f'(u)G'(v)^2)$.
We hence have $G''(v)v \in C^1(I)$, and (i) is verified.

From \eqref{dGfo}, we see that $G'(v)>0$ and
\begin{equation*}
h(v)=\frac{f'(u)F(u)}{f(u)^2}.
\end{equation*}
The equivalence of \eqref{aasfs} and \eqref{aashs}
then follows easily from this and \eqref{dG3fo}.
We also see from \eqref{dGfo} that $H(v)=F(u)-f(u)u$,
which gives
\begin{equation*}
H'(v)=\frac{du}{dv}\frac{d}{du}(F(u)-f(u)u)=G'(v)(f(u)u-f'(u)).
\end{equation*}
Hence \eqref{aasfw} and \eqref{aasGw} are equivalent.

It remains to prove (iii).
Taylor's theorem yields
$F(s)=\be_0^2+f'(\pm 1)(1-|s|)^2 +o((1-|s|)^2)$ as $s \to \pm 1$.
From this and the fact that $F(u)=v^2$,
we have
\begin{equation*}
\frac{\sqrt{\be_0-|v|}}{1-|u|}=\sqrt{\frac{\be_0^2-F(u)}{(\be_0+|v|)(1-|u|)^2}}
 \to \sqrt{\frac{-f'(\pm 1)}{2\be_0}} \quad (v \to \pm \be_0).
\end{equation*}
Using this and \eqref{dGfo} shows that as $v \to \pm \be_0$,
\begin{gather*}
\sqrt{\be_0-|v|} G'(v)
 =\frac{\sqrt{\be_0-|v|}}{1-|u|} \cdot \frac{1-|u|}{f(u)} \cdot v
  \to \sqrt{\frac{\be_0}{-2f'(\pm 1)}},
\\
(\be_0-|v|)^{3/2} G''(v) =\frac{(\be_0-|v|)^{3/2}}{(1-|u|)^3}
 \cdot \frac{(1-|u|)^3}{f(u)^3} \cdot (f(u)^2 -f'(u)F(u))
  \to \pm \sqrt{\frac{\be_0}{-8f'(\pm 1)}}.
\end{gather*}
Since $G'(v)>0$, we obtain \eqref{Gpas} and \eqref{Gdpas}.
Thus the lemma follows.
\end{proof}

\begin{proof}[Proof of Lemma~\ref{lem:bt1}]
We prove the lemma without using a spectral theory for compact operators 
or the Fredholm alternative.
We note that $T$ is symmetric,
since the right-hand side of \eqref{Tsym} vanish if $\varphi,\psi \in {\cX}$.

First we claim that the following holds:
\begin{equation}\label{phqdeq}
\mbox{if } \eta \in {\cY} \mbox{ and }  \varphi \in {\cZ}
 \mbox{ satisfy } \langle \varphi, T\psi \rangle=\langle \eta, \psi \rangle
  \mbox{ for any } \psi \in {\cX},
   \mbox{ then } \varphi \in {\cX} \mbox{ and } T \varphi=\eta.
\end{equation}
Suppose that $\eta \in {\cY}$, $\varphi \in {\cZ}$ and
$\langle \varphi, T\psi \rangle=\langle \eta, \psi \rangle$ for $\psi \in {\cX}$.
Then $\varphi$ satisfies the following equation in the distributional sense:
\begin{equation}\label{phqeq}
\varphi_{xx}+q(x) \varphi=\eta(x), \quad x \in (-1,1) \setminus \{0\}. 
\end{equation}
From a regularity theory for differential equations,
we see that $\varphi$ belongs to $X_0$ and satisfies \eqref{phqeq} in the classical sense.
Hence it follows from \eqref{Tsym} that for any $\psi \in {\cX}$,
\begin{align*}
0=&\psi(1) \varphi_x(1) -\psi(-1) \varphi_x(-1)
 -\psi_x(-0) \big\{ (\varphi (-0) +a\varphi_x (-0)) -(\varphi (+0) -a\varphi_x (+0))\big\}
\\
&+(\psi (-0) +a\psi_x (-0)) (\varphi_x(-0) -\varphi_x(+0)).
\end{align*}
Since $\psi \in {\cX}$ is arbitrary,
we find that $\varphi_x(1)=\varphi_x(-1)=0$,
$\varphi (-0) +a\varphi_x (-0)=\varphi (+0) -a\varphi_x (+0)$ and $\varphi_x(-0) =\varphi_x(+0)$.
Therefore $\varphi \in {\cX}$ and $T \varphi=\eta$, as claimed.

Next we show that
\begin{equation}\label{Ties}
\| \varphi \|_{\cZ} \le C\| T\varphi \|_{\cZ} 
\quad 
\mbox{for all } \varphi \in K(T)^\perp \cap {\cX},
\end{equation}
where $C>0$ is a constant.
Suppose that this claim were false. 
Then we could find a sequence $\{ \varphi_n \} \subset K(T)^\perp \cap {\cX}$ such that
$\| \varphi_n \|_{\cZ}=1$ and $T\varphi_n \to 0$ in ${\cZ}$ as $n \to \infty$.
Since $\| (\varphi_n)_{xx} \|_{\cZ} =\| T\varphi_n -q\varphi_n \|_{\cZ}
\le \| T\varphi_n \|_{\cZ} +\| q\|_{\cY} \| \varphi_n \|_{\cZ}$, 
we see that $\{ (\varphi_n)_{xx} \}$ is bounded in $\cZ$.
According to the Rellich-Kondrachov theorem,
there are a subsequence $\{ \tilde \varphi_n \}$ of $\{ \varphi_n \}$ 
and $\varphi \in {\cZ}$ such that $\tilde \varphi_n \to \varphi$ in ${\cZ}$.
Since $T$ is symmetric, we have
\begin{equation*}
\langle \varphi, T\psi \rangle =\lim_{n \to \infty} \langle \tilde \varphi_n, T\psi \rangle
 =\lim_{n \to \infty} \langle T\tilde \varphi_n, \psi \rangle =0
  \quad \mbox{for all } \psi \in {\cX}.
\end{equation*}
By \eqref{phqdeq}, we infer that $\varphi \in {\cX}$ and $T\varphi=0$.
Hence $\varphi \in K(T)$.
On the other hand, 
since $\varphi_n$ satisfies $\| \varphi_n \|_{\cZ}=1$ and $\varphi_n \in K(T)^\perp$,
we have $\| \varphi \|_{\cZ}=1$ and $\varphi \in K(T)^\perp$.
This leads to a contradiction, because no nonzero function exists in $K(T) \cap K(T)^\perp$.
Therefore \eqref{Ties} holds.

Finally let us prove \eqref{ReKp}.
We have $R(T)^\perp=K(T)$,
since using \eqref{phqdeq} for $\eta=0$ yields $R(T)^\perp \subset K(T)$
and the fact that $T$ is symmetric gives $K(T) \subset R(T)^\perp$.
From this it follows that $\overline{R(T)}=K(T)^\perp$,
where $\overline{R(T)}$ is the closure of $R(T)$ in $\cZ$.
Therefore we only need to show that
\begin{equation}\label{oRYR}
\overline{R(T)} \cap {\cY} =R(T).
\end{equation}
Since it is clear that $R(T) \subset \overline{R(T)} \cap {\cY}$,
what is left is to prove $\overline{R(T)} \cap {\cY} \subset R(T)$.
Let $\eta \in \overline{R(T)} \cap {\cY}$.
By definition, there exists a sequence $\{ \varphi_n \} \subset {\cX}$ 
such that $T\varphi_n \to \eta$ in $\cZ$ as $n \to \infty$.
Since ${\cX}=(K(T)^\perp \cap {\cX}) \oplus K(T)$,
we can take $\tilde \varphi_n\in K(T)^\perp \cap {\cX}$ 
such that $\varphi_n -\tilde \varphi_n \in K(T)$.
Then 
\begin{equation}\label{typl}
T\tilde \varphi_n=T\varphi_n \to \eta \quad \mbox{in } \cZ \mbox{ as } n \to \infty.
\end{equation}
This together with \eqref{Ties} shows that $\{ \tilde \varphi_n \}$ is a Cauchy sequence in $\cZ$.
Hence $\tilde \varphi_n \to \varphi$ in $\cZ$ for some $\varphi \in {\cZ}$.
Combining this with \eqref{typl} and the fact that $T$ is symmetric,
we deduce that
\begin{equation*}
\langle \varphi, T\psi \rangle =\lim_{n \to \infty} \langle \tilde \varphi_n, T\psi \rangle
 =\lim_{n \to \infty} \langle T\tilde \varphi_n, \psi \rangle =\langle \eta, \psi \rangle
  \quad \mbox{for all } \psi \in {\cX}.
\end{equation*}
We see from \eqref{phqdeq} that $\varphi \in {\cX}$ and $T\varphi=\eta$,
which means $\eta \in R(T)$. 
Thus \eqref{oRYR} holds, and the proof is complete.
\end{proof}



\begin{thebibliography}{99}
\bibitem{AN13} 
R.~Adami and D.~Noja,
{\it Stability and symmetry-breaking bifurcation for the ground states of a NLS 
with a $\del'$ interaction}, 
Comm. Math. Phys. {\bf 318} (2013), no.~1, 247--289.

\bibitem{AN14} 
R.~Adami and D.~Noja,
{\it Exactly solvable models and bifurcations: the case of the cubic NLS 
with a $\del$ or a $\del'$ interaction in dimension one},
Math. Model. Nat. Phenom. {\bf 9} (2014), no.~5, 1--16.

\bibitem{AG18} 
J.~Angulo Pava and N.~Goloshchapova, 
{\it Extension theory approach in the stability of the standing waves 
for the NLS equation with point interactions on a star graph},
Adv. Differential Equations {\bf 23} (2018),  no.~11-12, 793--846.

\bibitem{AG20} 
J.~Angulo Pava and N.~Goloshchapova, 
{\it Stability properties of standing waves for NLS equations with the $\del'$-interaction}, 
Phys. D {\bf 403} (2020), 132332, 24 pp.

\bibitem{A64}
F.~V.~Atkinson, 
{\it Discrete and continuous boundary problems},
Mathematics in Science and Engineering, Vol. 8, Academic Press, New York-London, 1964.

\bibitem{CH78}
R.~G.~Casten and C.~J.~Holland,
{\it Instability results for reaction diffusion equations with Neumann boundary conditions},
J. Differential Equations {\bf 27}  (1978), no.~2, 266--273.

\bibitem{CR71}
M.~G.~Crandall and P.~H.~Rabinowitz, 
{\it Bifurcation from simple eigenvalues},
J. Functional Analysis {\bf 8} (1971), 321--340.

\bibitem{F90} 
Q.~Fang,
{\it Asymptotic behavior and domain-dependency of solutions
to a class of reaction-diffusion systems with large diffusion coefficients},
Hiroshima Math. J. {\bf 20} (1990), no.~3, 549--571.

\bibitem{GO20} 
N.~Goloshchapova and M.~Ohta,
{\it Blow-up and strong instability of standing waves 
for the NLS-$\del$ equation on a star graph},
Nonlinear Anal. {\bf 196} (2020), 111753, 23 pp.

\bibitem{HR85} 
J.~K.~Hale and C.~Rocha, 
{\it Bifurcations in a parabolic equation with variable diffusion},
Nonlinear Anal.  {\bf 9}  (1985),  no.~5, 479--494.

\bibitem{HV84} 
J.~K.~Hale and J.~Vegas,
{\it A nonlinear parabolic equation with varying domain},
Arch. Rational Mech. Anal. {\bf 86} (1984), no.~2, 99--123.

\bibitem{KST18}
R.~Kajikiya, I.~Sim and S.~Tanaka,
{\it Symmetry-breaking bifurcation for the Moore-Nehari differential equation},
NoDEA Nonlinear Differential Equations Appl. {\bf 25} (2018), no.~6, Paper No.~54, 22 pp.

\bibitem{K}
T.~Kan,
{\it Semilinear elliptic equations on thin dumbbell-shaped domains},
in preparation. 


\bibitem{K12}
H.~Kielh\"{o}fer,
{\it Bifurcation theory: 
An introduction with applications to partial differential equations
2nd edition},
Applied Mathematical Sciences 156, Springer, New York.

\bibitem{Ma79}
H.~Matano,
{\it Asymptotic behavior and stability of solutions of semilinear diffusion equations},
Publ. Res. Inst. Math. Sci.  {\bf 15}  (1979), no.~2, 401--454.

\bibitem{MEF91} 
M.~Mimura, S.~Ei and Q.~Fang,
{\it Effect of domain-shape on coexistence problems in a competition-diffusion system},
J. Math. Biol. {\bf 29} (1991), no.~3, 219--237.

\bibitem{Mo90} 
Y.~Morita,
{\it Reaction-diffusion systems in nonconvex domains: invariant manifold and reduced form},
J. Dynam. Differential Equations {\bf 2} (1990), no.~1, 69--115.

\bibitem{MJ92} 
Y.~Morita and S.~Jimbo,
{\it Ordinary differential equations (ODEs) on inertial manifolds
for reaction-diffusion systems in a singularly perturbed domain with several thin channels},
J. Dynam. Differential Equations {\bf 4} (1992), no.~1, 65--93.

\bibitem{O61} 
Z.~Opial, 
{\it Sur les p\'{e}riodes des solutions de l'\'{e}quation diff\'{e}rentielle $x''+g(x)=0$},
Ann. Polon. Math. {\bf 10} (1961), 49--72.

\bibitem{S90}
R.~Schaaf,
{\it Global Solution Branches of Two Point Boundary Value Problems},
Lecture Notes in Mathematics, Vol. 1458, Springer-Verlag, Berlin,  1990.

\bibitem{ST19}
I.~Sim and S.~Tanaka,
{\it Symmetry-breaking bifurcation for the one-dimensional H\'{e}non equation},
Commun. Contemp. Math.  {\bf 21}  (2019),  no.~1, 1750097, 24 pp.

\bibitem{SW81} 
J.~Smoller and A.~Wasserman, 
{\it Global bifurcation of steady-state solutions},
J. Differential Equations {\it 39} (1981), no.~2, 269--290.

\bibitem{T13}
S.~Tanaka,
{\it Morse index and symmetry-breaking for positive solutions of 
one-dimensional H\'{e}non type equations},
J. Differential Equations  {\bf 255}  (2013),  no.~7, 1709--1733.

\bibitem{T17}
S.~Tanaka,
{\it Symmetry-breaking bifurcation for the one-dimensional Liouville type equation}, 
J. Differential Equations {\bf 263} (2017), no.~10, 6953--6973.

\bibitem{V83} 
J.~Vegas,
{\it Bifurcations caused by perturbing the domain in an elliptic equation},
J. Differential Equations {\bf 48} (1983), no.~2, 189--226.

\bibitem{W83} 
J.~Waldvogel, 
{\it The Period in the Volterra-Lotka Predator-Prey Model},
SIAM J. Numer. Anal. {\bf 20} (1983), no.~6, 1264--1272.

\bibitem{W86} 
J.~Waldvogel, 
{\it The period in the Lotka-Volterra system is monotonic},
J. Math. Anal. Appl. {\bf 114} (1986), no.~1, 178--184.
\end{thebibliography}
\end{document}